# Ordered groups of formal series, and a conjugacy problem

by Vincent Bagayoko

imj-prg

*Email:* bagayoko@imj-prg.fr

**Abstract**

Given an ordered field $\mathbb{T}$ of formal series over an ordered field $\mathbf{R}$ equipped with a composition law $\circ \colon \mathbb{T} \times \mathbb{T}^{>\mathbf{R}} \longrightarrow \mathbb{T}$, we give conditions for $(\mathbb{T}^{>\mathbf{R}}, \circ)$ to be a group. We show that classical fields of transseries satisfy these conditions. We then give further conditions on $\mathbb{T}$ under which $(\mathbb{T}^{>\mathbf{R}}, \circ, <)$ is a linearly ordered group with exactly three conjugacy classes and solve the open problem of existence of such a linearly ordered group.

# Introduction

**The conjugacy problem for linearly ordered groups** We call *conjugacy problem* the existence of a linearly (bi-)ordered, non-trivial group $(\mathcal{G}, \cdot, 1, <)$ in which any two strictly positive elements are conjugate, i.e. such that

$$\forall f, g \in \mathcal{G} \, (f, g > 1 \Longrightarrow \exists h \in \mathcal{G} \, (h\, f = g\, h)). \tag{1}$$

To the best of our knowledge, the existence of such a linearly bi-ordered group is an open problem (see [14, Problem 3.31]). In this paper, we propose a solution.

A straightforward approach to that problem is to start with an ordered group $\mathcal{G}$, to fix a positive element $f > 1$ in $\mathcal{G}$, and for each $g \in \mathcal{G}$ which is not conjugate of $f$, to choose an ordering on the HNN extension of $\mathcal{G}$ conjugating $g$ with $f$, which extends the ordering on $\mathcal{G}$. Then proceed inductively until all conjugacy equations are solved. Unfortunately, there may exist several such extensions, or no such extensions at all [13]. A sounder approach is to restrict the class of ordered groups in which one is looking for a solution. It is then a matter of finding groups with sufficiently rigid properties that ordering words involving conjugating elements $h$ in (1) is easy.

Our solution is a group under composition of generalised power series called *hyperseries*. Generalised power series [22] are formal sums

$$s = \sum_{\mathfrak{m} \in \mathfrak{M}} s(\mathfrak{m}) \, \mathfrak{m},$$

with coefficients $s(\mathfrak{m})$ in an ordered field $\mathbf{R}$ and monomials $\mathfrak{m}$ in a linearly ordered Abelian group $(\mathfrak{M}, \cdot, 1, \prec)$, whose support $\{\mathfrak{m} \in \mathfrak{M} : s(\mathfrak{m}) \neq 0\}$ is well-ordered in $(\mathfrak{M}, \succ)$. The set $\mathbf{R}[\![\mathfrak{M}]\!]$ of all such series is an ordered valued field in a natural way (see Section 1.1). A very useful feature of $\mathbf{R}[\![\mathfrak{M}]\!]$ is a form of formal completeness, such that it can be construed as a formal analog of a Banach algebra [24]. Finding approximate solutions to equations over $\mathbf{R}[\![\mathfrak{M}]\!]$ often brings one half of the way toward obtaining an exact solution of that equation. It is tempting to apply this to conjugacy equations.

**Constructing groups of series under composition**  The ordered field $\mathbb{S} = \mathbb{R}[\![x^{\mathbb{Z}}]\!]$ of formal Laurent series over $\mathbb{R}$, where $x > \mathbb{R}$, has a composition law

$$\left(\sum_{k \leqslant n_0} r_k\, x^k\right) \circ \left(\sum_{m \leqslant n_1} s_m\, x^m\right) = \sum_{k \leqslant n_0} r_k \left(\sum_{m_1, \ldots, m_k \leqslant n_1} s_{m_1} \cdots s_{m_k}\right) x^{m_1 + \cdots + m_k}.$$

which is defined if $\sum_{m \leqslant n_1} s_m\, x^m$ is positive infinite. It can be shown that each $\hat{f}$ for $f \in \mathbb{S}^{>\mathbb{R}}$ is strictly increasing, so $(\mathbb{S}^{>\mathbb{R}}, \circ, x, <)$ is a linearly ordered monoid. But $\mathbb{S}^{>\mathbb{R}}$ is not a group, for $x^2$ has no compositional inverse. An inverse exists in the extension $\mathbb{R}[\![x^{\mathbb{Q}}]\!]$ of $\mathbb{S}$, as the series $\sqrt{x} = x^{1/2}$, and $(\mathbb{R}[\![x^{\mathbb{Q}}]\!]^{>\mathbb{R}}, \circ, x, <)$ is an ordered group for the natural extension of $\circ$. However, this is not a solution to the conjugacy problem, since $x^2$ and $x+1$ are both positive (i.e. larger than $x$) but not conjugate. A natural conjugating element exists in extensions of $\mathbb{R}[\![x^{\mathbb{Q}}]\!]$ called fields of transseries [19, 23, 30], i.e. series involving combinations of exponentials and logarithms. There is a field $\mathbb{R}\langle\!\langle \omega \rangle\!\rangle$ of so-called $\omega$-series which is closed under derivation and composition with positive infinite elements on the right. It was constructed independently by Schmeling [30] as a field of transseries, and by Berarducci and Mantova [11, 12] as a subfield of Conway's class **No** of surreal numbers [15]. In $\mathbb{R}\langle\!\langle \omega \rangle\!\rangle$, the smallest infinite ordinal $\omega \in \mathbf{On} \subseteq \mathbf{No}$ plays the role of the positive infinite variable $x$. The $\omega$-series $\omega^2$ and $\omega + 1$ are conjugate, via

$$\left(\frac{1}{\log 2} \log \log \omega\right) \circ \omega^2 = \omega \circ \left(\frac{1}{\log 2} \log \log \omega\right) + 1.$$

However, there is no $\omega$-series solution $f$ of Abel's conjugacy equation

$$f \circ (\omega + 1) = \exp(\omega) \circ f. \tag{2}$$

It is possible to add formal solutions to such equations while preserving all properties of $\omega$-series. The abstract notion of transserial field (see Section 2) is malleable enough to encompass extensions of $\mathbb{R}\langle\!\langle \omega \rangle\!\rangle$ with solutions to (2) and more general conjugacy equations. Consider the transserial field $\mathbb{L}$ of logarithmic hyperseries [17]. Those are series over $\mathbb{R}$ whose building blocks are symbols $\ell_\gamma$ for arbitrary ordinals $\gamma \in \mathbf{On}$, such that each $\ell_{\omega^{\gamma+1}}$ for is a privileged solution to the conjugacy equation

$$\ell_{\omega^{\gamma+1}} \circ \ell_{\omega^\gamma} = \ell_{\omega^\gamma} - 1$$

between the negative elements $\ell_{\omega^\gamma}$ and $\ell_0 - 1$ of the ordered monoid $(\mathbb{L}^{>\mathbb{R}}, \circ, \ell_0, <)$. We showed [9] that there is a minimal extension $\tilde{\mathbb{L}}$ of $\mathbb{L}$, called the field of *finitely nested hyperseries*, which contains compositional inverses $e_\gamma^{\ell_0}$ of each $\ell_\gamma$ for $\gamma \in \mathbf{On}$. We then extended [5] the derivation and composition law on $\mathbb{L}$ to $\tilde{\mathbb{L}}$.

The subclass $\tilde{\mathbb{L}}^{>\mathbb{R}}$ of $\tilde{\mathbb{L}}$ of positive infinite elements constitutes the first candidate of a solution to the conjugacy problem. It was unknown yet whether $(\tilde{\mathbb{L}}^{>\mathbb{R}}, \circ, \ell_0)$ and other smaller sets of transseries form groups. More globally, a general result for showing that a transserial field $\mathbb{S}$ with a composition law $\circ$ induces a group $(\mathbb{S}^{>\mathbf{R}}, \circ, x)$ is lacking. We provide one in this paper. We obtain as corollaries that:

**Theorem 1.** [Corollary 3.18] *The monoid $(\tilde{\mathbb{L}}^{>\mathbb{R}}, \circ, \ell_0)$ is a group.*

**Theorem 2.** [Corollary 3.20] *The monoid $(\mathbb{R}\langle\!\langle \omega \rangle\!\rangle^{>\mathbb{R}}, \circ, \omega)$ is a group.*

We then propose a method for solving conjugacy equations in such groups and apply it to show that $(\tilde{\mathbb{L}}^{>\mathbb{R}}, \circ, \ell_0)$ solves the conjugacy problem:

**Theorem 3.** [Corollary 4.22] *Any two positive elements in the ordered group $(\tilde{\mathbb{L}}^{>\mathbb{R}}, \circ, \ell_0, <)$ are conjugate.*

Our results will apply to **No** provided a suitable derivation and composition law are defined on this field (see [6, Conclusion]).

**Outline of the paper** We start in Section 1 by introducing our notion of field of series. We also gather basic tools for solving equations and defining functions on those fields. We recall well-known results regarding the formal notion of *summability* (see Section 1.2).

In Section 2, we define *transserial fields* (Section 2.1), which include the fields of $\omega$-series, finitely nested hyperseries, and surreal numbers. We introduce the operations of exponentiation (Section 2.2) and powering (Section 2.3) and the ordering of steepness (Section 2.4) on such fields. Given a transserial field $\mathbb{T} = \mathbf{R}[\![\mathfrak{M}]\!]$, we specify our constraints on derivations $\partial$ and composition laws $\circ$ on $\mathbb{T}$ (Sections 2.5 and 2.6). We also rely (Section 2.7) on a notion of *Taylor expansion* that appears, more or less explicitly, in all previous work related to transseries [19, 23, 30, 18, 25].

In Section 3, we tackle our first goal of constructing compositional inverses in transserial fields $\mathbb{T} = \mathbf{R}[\![\mathfrak{M}]\!]$ equipped with a derivation $\partial$ and composition law $\circ$ with identity $x$. Our main result Theorem 3.1 consists in giving conditions under which the class

$$\mathcal{G} := \{x + \delta : r\,\delta < x \text{ for all } r \in \mathbf{R}\}$$

is a group under composition. In order to find the inverse of a series $x + \delta \in \mathcal{G}$, we represent it (Section 3.4) as a possibly transfinite composition

$$x + \delta = (x + \delta_0) \circ (x + \delta_1) \circ \cdots \circ (x + \delta_\gamma) \circ \cdots$$

where $\gamma$ ranges in an ordinal $\lambda$, and the sequence $(\delta_\gamma)_{\gamma < \lambda}$ is strictly increasing for the steepness ordering. We obtain (Section 3.5) the inverse $(x + \delta)^{\mathrm{inv}}$ of $x + \delta$ as an infinite composition in the other direction

$$\cdots \circ (x + \delta_\gamma)^{\mathrm{inv}} \circ \cdots \circ (x + \delta_1)^{\mathrm{inv}} \circ (x + \delta_0)^{\mathrm{inv}} = (x + \delta)^{\mathrm{inv}},$$

where each $(x + \delta_\gamma)^{\mathrm{inv}}$ can be computed directly using Taylor expansions (see Section 3.3). In certain cases, one deduces from the fact that $\mathcal{G}$ is a group that $(\mathbb{T}^{>\mathbf{R}}, \circ, x)$ is a group, yielding Theorems 1 and 2.

Section 4 focuses on the conjugacy problem. We generalise Edgar's method [21] for solving conjugacy equations over log-exp transseries, as is shown in Sections 4.1 and 4.2, so that it apply to the case of finitely nested hyperseries, as shown in Section 4.3. Our general theorem Theorem 4.6 implies Theorem 3, and it will also apply in a straightforward way to the case of surreal numbers.

Sections 4.4 and 4.5 are dedicated to studying the simple inequality

$$f \circ g \geqslant g \circ f \qquad (3)$$

for positive elements in the ordered group $(\mathbb{T}^{>\mathbf{R}}, \circ, x, <)$, provided this group solves the conjugacy problem. This very elementary inequality gives a lot of information on that ordered group. Solving (3) entails determining centralizers $\mathcal{C}(f) = \{h \in \mathbb{T}^{>\mathbf{R}} : h \circ f = f \circ h\}$ of elements $f \in \mathbb{T}^{>\mathbf{R}}$, which can be done by studying iterates $f^{[r]}$ of $f$, where $r$ ranges in $\mathbf{R}$. For fixed $r \in \mathbf{R}$, the function $f \mapsto f^{[r]}$ is a quasi-endormorphism of $\mathbb{T}^{>\mathbf{R}}$ in the sense of [28] (see Propositions 4.27 and 4.28). More importantly, for fixed $f \neq x$, the function $r \mapsto f^{[r]}$ is an isomorphism $(\mathbf{R}, +, 0, <) \longrightarrow (\mathcal{C}(f), \circ, x, <)$, thus

**Proposition 4.** [Proposition 4.27] *We have $(\mathcal{C}(f), \circ, x, <) \simeq (\mathbf{R}, +, 0, <)$ for all $f \in \mathbb{T}^{>\mathbf{R}}$ with $f \neq x$.*

Using easy computations with iterates, we derive important first-order properties of $(\mathbb{T}^{>\mathbf{R}}, \circ, x, <)$. In particular, we show:

**Proposition 5.** [Propositions 4.29 and 4.30] *For all $f, g > x$, we have*

$$f > \mathcal{C}(g) \Longrightarrow f \circ g > g \circ f \qquad \text{and} \qquad f \geqslant g \Longrightarrow \forall g_0 \in \mathcal{C}(g)\, (\exists f_0 \in \mathcal{C}(f)\, (f_0 \geqslant g_0)).$$

Those last results imply that $(\mathbb{T}^{>\mathbf{R}}, \circ, x, <)$, alongside certain groups of germs definable in o-minimal structures, is a growth order group in the sense of [8].

**Conventions.** Before we start, we set a few conventions.

**Set theory.** We adopt the set-theoretic framework of [9]. In particular, the underlying set theory of this paper is NBG set theory. This is a conservative extension of ZFC which allows us to prove statements about proper classes. We will use boldface font for classes that may not be sets, and reserve the standard font for sets (with the notable exceptions of the fraktur and mathbb fonts $\mathfrak{M}, \mathfrak{N}, \mathbb{T}, \mathbb{S} \ldots$).

**Ordered monoids.** If $(\mathbf{M}, +, 0, <)$ is an ordered monoid such as $\mathbb{N}$, $\mathbb{R}$, one of our hyperserial fields $\mathbb{T}$ or groups of monomials $\mathfrak{M}$, then $\mathbf{M}^>$ denotes its subclass of strictly positive elements in $\mathbf{M}$, whereas $\mathbf{M}^{\neq}$ denotes the class of non-zero elements.

# 1 Strongly linear algebra

Here, we recall well-known notions on fields of Hahn series. See [24, 6, 10] for more detail.

## 1.1 Fields of well-based series

Let $\mathbf{R}$ be an ordered field. Let $(\mathfrak{M}, \cdot, 1, \prec)$ be a linearly ordered Abelian group. We write $\mathbf{R}[\![\mathfrak{M}]\!]$ for the class of functions $s : \mathfrak{M} \longrightarrow \mathbf{R}$ whose support

$$\operatorname{supp} s := \{\mathfrak{m} \in \mathfrak{M} : s(\mathfrak{m}) \neq 0\}$$

is a *well-based* set, i.e. a set which is well-ordered in the reverse order $(\mathfrak{M}, \succ)$. We allow $\mathfrak{M}$ to be a proper class, but supports of series in $\mathbb{S}$ should always be sets. We see elements $s$ of $\mathbb{S}$ as formal *well-based series*

$$s = \sum_{\mathfrak{m} \in \mathfrak{M}} s(\mathfrak{m})\, \mathfrak{m}. \tag{1.1}$$

By [22], the class $\mathbb{S}$ is a field under the pointwise sum

$$(s + t) := \sum_{\mathfrak{m}} (s(\mathfrak{m}) + t(\mathfrak{m}))\, \mathfrak{m},$$

and the Cauchy product

$$s\, t := \sum_{\mathfrak{m}} \left( \sum_{\mathfrak{u}\mathfrak{v} = \mathfrak{m}} s(\mathfrak{u})\, t(\mathfrak{v}) \right) \mathfrak{m},$$

(where each sum $\sum_{\mathfrak{u}\mathfrak{v}=\mathfrak{m}} s(\mathfrak{u})\, t(\mathfrak{v})$ has finite support). We have an ordering $<$ on $\mathbb{S}$ whose positive cone $\mathbb{S}^> = \{s \in \mathbb{S} : s > 0\}$ is defined by

$$\mathbb{S}^> := \{s \in \mathbb{S} : s \neq 0 \land s(\max \operatorname{supp} s) > 0\}.$$

We say that $\mathbb{S}$ is an *ordered field of well-based series*. We set $|s| := \max(s, -s)$ for all $s \in \mathbb{S}$. We have an embedding of ordered groups

$$(\mathfrak{M}, \cdot, 1, \prec) \longrightarrow (\mathbb{S}^{>}, \cdot, 1, <)\,;\, \mathfrak{m} \mapsto \sum_{\mathfrak{n} = \mathfrak{m}} \mathfrak{n}$$

and we identify $\mathfrak{M}$ with its image in $\mathbb{S}^{>}$. The elements of $\mathfrak{M}$ are called *monomials*, whereas those in $\mathbf{R}^{\times}\mathfrak{M}$ are called *terms*. We write

$$\operatorname{term} s := \{s(\mathfrak{m})\,\mathfrak{m} : \mathfrak{m} \in \operatorname{supp} s\}$$

for the set of terms of $s$. If $\operatorname{supp} s \neq \varnothing$, then we write

$$\begin{aligned}\mathfrak{d}(s) &:= \max \operatorname{supp} s \in \mathfrak{M} \quad \text{and} \\ \tau_s &:= s(\mathfrak{d}(s))\,\mathfrak{d}(s) \in \mathbf{R}^{\times}\mathfrak{M}\end{aligned}$$

respectively for the *dominant monomial* and *dominant term* of $s$. For $s, t \in \mathbb{S}$, we say that $t$ is a *truncation* of $s$ and we write $t \trianglelefteq s$ if $\operatorname{supp}(s - t) \succ \operatorname{supp} s$ and $s \triangleleft t$ if $s \neq t$ and $s \trianglelefteq t$.

The ordering on $\mathfrak{M}$ extends into a strict ordering $\prec$ on $\mathbb{S}$ defined by $s \prec t$ if and only if $\mathbf{R}^{>}|s| < |t|$. We write $s \preccurlyeq t$ if $t \prec s$ is false, i.e. if there is $r \in \mathbf{R}^{>}$ with $|s| \leqslant r\,|t|$. We also write $s \asymp t$ if $s \preccurlyeq t$ and $t \preccurlyeq s$, i.e. if there is an $r \in \mathbf{R}^{>}$ with $r\,|s| \geqslant |t|$ and $r\,|t| \geqslant s$. Finally, we write $s \sim t$ if $s - t \prec s$, or equivalently $s - t \prec t$. The relation $\sim$ is an equivalence relation on $\mathbb{S}^{\times}$. The relation $\preccurlyeq$ is a dominance relation as per [4, Definition 3.1.1], and it corresponds to the natural valuation on the ordered field $(\mathbb{S}, +, \times, <)$. In particular $(\mathbb{S}, +, \cdot, <, \prec)$ is an ordered valued field with convex valuation ring $\mathbb{S}^{\preccurlyeq} := \{s \in \mathbb{S} : s \preccurlyeq 1\}$. When $s, t$ are non-zero, we have $s \prec t$ (resp. $s \preccurlyeq t$, resp. $s \asymp t$) if and only if $\mathfrak{d}(s) \prec \mathfrak{d}(t)$ (resp. $\mathfrak{d}(s) \preccurlyeq \mathfrak{d}(t)$, resp. $\mathfrak{d}(s) = \mathfrak{d}(t)$).

If $\mathfrak{G} \subseteq \mathfrak{M}$ is a subclass, then we write $\mathbf{R}[\![\mathfrak{G}]\!]$ for the class of elements $s \in \mathbb{S}$ with $\operatorname{supp} s \subseteq \mathfrak{G}$. Note that this is a subgroup of $\mathbb{S}$, and that it is a subfield if and only if $\mathfrak{G}$ is a group. We write

$$\begin{aligned}\mathbb{S}_{\succ} &:= \mathbf{R}[\![\mathfrak{M}^{\succ}]\!] = \{s \in \mathbb{S} : \operatorname{supp} s \succ 1\} \\ \mathbb{S}^{\prec} &:= \mathbf{R}[\![\mathfrak{M}^{\prec}]\!] = \{s \in \mathbb{S} : s \prec 1\}, \quad \text{and} \\ \mathbb{S}^{>,\succ} &:= \{s \in \mathbb{S} : s > \mathbf{R}\} = \{s \in \mathbb{S} : s \geqslant 0 \wedge s \succ 1\}.\end{aligned}$$

Series in $\mathbb{S}^{\prec}$ and $\mathbb{S}^{>,\succ}$ are respectively said *infinitesimal*, and *positive infinite*.

## 1.2 Summable families

Let $\mathbf{I}$ be a class. A family $(s_i)_{i \in I}$ in $\mathbb{S}$ is said *summable* if

i. $\bigcup_{i \in \mathbf{I}} \operatorname{supp} s_i$ is a well-based set, and

ii. the set $\{i \in \mathbf{I} : \mathfrak{m} \in \operatorname{supp} s_i\}$ is finite for all $\mathfrak{m} \in \mathfrak{M}$.

Then we may define the sum $\sum_{i \in \mathbf{I}} s_i$ of $(s_i)_{i \in \mathbf{I}}$ as the series

$$\sum_{i \in \mathbf{I}} s_i := \sum_{\mathfrak{m}} \left( \sum_{i \in \mathbf{I}} (s_i)(\mathfrak{m}) \right) \mathfrak{m}.$$

**Example 1.1.** Given $s \in \mathbb{S}$, the family $(s(\mathfrak{m})\,\mathfrak{m})_{\mathfrak{m} \in \mathfrak{M}}$ is summable with sum $s$. This lends a formal meaning to the identification in (1.1).

**Example 1.2.** *As a consequence of Neumann's theorems [29], for all $(r_k)_{k \in \mathbb{N}} \in \mathbf{R}^{\mathbb{N}}$ and $\varepsilon \prec 1$ the family $(r_k \varepsilon^k)_{k \in \mathbb{N}}$ is summable. This allows for evaluating formal power series [10] and real-analytic functions [16] (when $\mathbb{R} = \mathbf{R}$) on subclasses of $\mathbb{S}$.*

### 1.3 Strong linearity

Let $\mathbb{S} = \mathbf{R}[\![\mathfrak{M}]\!]$ and $\mathbb{T} = \mathbf{R}[\![\mathfrak{N}]\!]$ be ordered fields of well-based series. Consider an $\mathbf{R}$-linear function $\Phi : \mathbb{S} \longrightarrow \mathbb{T}$. Then $\Phi$ is *strongly linear* if for every summable family $(s_i)_{i \in I}$ in $\mathbb{S}$, the family $(\Phi(s_i))_{i \in I}$ in $\mathbb{T}$ is summable, with $\Phi(\sum_{i \in I} s_i) = \sum_{i \in I} \Phi(s_i)$. By [24, Proposition 3.5], the function $\Phi$ is strongly linear if and only if for each $s \in \mathbb{S}$, the family $(\Phi(\mathfrak{m}))_{\mathfrak{m} \in \mathrm{supp}\, s}$ is summable with $\Phi(s) = \sum_{\mathfrak{m} \in \mathrm{supp}\, s} s_{\mathfrak{m}} \Phi(\mathfrak{m})$.

**Notation 1.3.** *Given a function $\Phi : \mathbf{X} \longrightarrow \mathbf{X}$ on a class $\mathbf{X}$ and a $k \in \mathbb{N}$, we will frequently write $\Phi^{[k]}$ for the k-fold iterate of $\Phi$. So $\Phi^{[k]}$ is the function $\mathbf{X} \longrightarrow \mathbf{X}$ with $\Phi^{[0]} = \Phi$ and $\Phi^{[k+1]} := \Phi^{[k]} \circ \Phi = \Phi \circ \Phi^{[k]}$ for all $k \in \mathbb{N}$.*

**Proposition 1.4.** *(corollary of [24, Theorem 6.2]) Let $\phi : \mathbb{S} \longrightarrow \mathbb{S}$ be strongly linear with $\phi(\mathfrak{m}) \prec \mathfrak{m}$ for all $\mathfrak{m} \in \mathfrak{M}$. Let $(r_k)_{k \in \mathbb{N}} \in \mathbf{R}^{\mathbb{N}}$. Then for all $s \in \mathbb{S}$, the family $(r_k \phi^{[k]}(s))_{k \in \mathbb{N}}$ is summable, and the function*

$$\sum_{k \in \mathbb{N}} r_k \phi^{[k]} : \mathbb{S} \longrightarrow \mathbb{S} ; s \mapsto \sum_{k \in \mathbb{N}} r_k \phi^{[k]}(s)$$

*is strongly linear.*

**Corollary 1.5.** *[24, Corollary 1.4] Let $\phi$ be as in Proposition 1.4. The function*

$$\mathrm{Id}_{\mathbb{S}} + \phi : \mathbb{S} \longrightarrow \mathbb{S} ; s \mapsto s + \phi(s)$$

*is bijective, with functional inverse $\sum_{k \in \mathbb{N}} (-1)^k \phi^{[k]}$.*

## 2 Transserial fields

We now introduce the main class of objects to which our results apply. The notion of transserial field used here is similar to that of transseries field in [30, Definition 2.2.1] where the difference is that we do not ask that Schmeling's axiom **T4** be satisfied. It is also close to the notion of pre-logarithmic Hahn field of [26, Definition 2.7].

### 2.1 Transserial fields

Let $\mathbf{R}$ be an ordered field and let $(\mathfrak{N}, \cdot, 1, \prec)$ be a linearly ordered Abelian group. We set $\mathbb{T} := \mathbf{R}[\![\mathfrak{N}]\!]$.

**Definition 2.1.** *Let $\log : (\mathbb{T}^>, \cdot, 1, <) \longrightarrow (\mathbb{T}, +, 0, <)$ be an embedding of ordered groups. We say that $(\mathbb{T}, \log)$ is a **transserial field** if*

- *$\log \mathfrak{N} \subseteq \mathbb{T}_{\succ}$,*
- *$\log s \leqslant s - 1$ for all $s \in \mathbb{T}^>$,*
- *$\log(1 + \varepsilon) = \sum_{k > 0} \frac{(-1)^{k+1}}{k} \varepsilon^k$ for all $\varepsilon \prec 1$,*
- *$\log \upharpoonright \mathbf{R}^> : (\mathbf{R}^>, \cdot, 1, <) \longrightarrow (\mathbf{R}, +, 0, <)$ is an isomorphism of ordered groups.*

*We call log the *logarithm* on $\mathbb{T}$.*

**Example 2.2.** Suppose that $\log_{\mathbf{R}} : \mathbf{R}^{>} \longrightarrow \mathbf{R}$ is an isomorphism of ordered groups with $\log_{\mathbf{R}}(r) \leqslant r + 1$ for all $r \in \mathbf{R}^{>}$. Let $\alpha$ be a non-zero limit ordinal. A small non-trivial example of transserial field over $\mathbf{R}$ is the field $\mathbb{T}_{\mathbf{R}, <\alpha} := \mathbf{R}[\![\mathfrak{L}_{\mathbf{R}, <\alpha}]\!]$ where $\mathfrak{L}_{\mathbf{R}, <\alpha}$ is a multiplicative copy of the group $\mathbf{R}^{\alpha} = \bigotimes_{\gamma < \alpha} (\mathbf{R}, +, 0)$ ordered lexicographically with prevalence on the first coordinates. Elements in $\mathfrak{L}_{\mathbf{R}, <\alpha}$ can be represented as formal products

$$\mathfrak{l} = \prod_{\gamma < \alpha} \ell_{\gamma}^{\mathfrak{l}_{\gamma}},$$

where $(\mathfrak{l}_{\gamma})_{\gamma < \alpha} \in \mathbf{R}^{\alpha}$ is a sequence and each $\ell_{\gamma}, \gamma < \alpha$ is a formal symbol. The group operation is the pointwise product $\mathfrak{l}\mathfrak{j} := \prod_{\gamma < \alpha} \ell_{\gamma}^{\mathfrak{l}_{\gamma} + \mathfrak{j}_{\gamma}}$. The lexicographic ordering has $\{\mathfrak{l} \in \mathfrak{L}_{\mathbf{R}, <\alpha} : \mathfrak{l}_{\eta} > 0 \text{ where } \eta = \min\{\gamma \in \mathbb{N} : \mathfrak{l}_{\gamma} \in \mathbf{R}^{\times}\}\}$ as a cone of strictly positive elements. Set $\log(\mathfrak{l}) := \sum_{\gamma < \alpha} \mathfrak{l}_{\gamma} \ell_{\gamma+1} \in \mathbb{T}_{\mathbf{R}, <\alpha}$ for all $\mathfrak{l} = \prod_{\gamma < \alpha} \ell_{\gamma}^{\mathfrak{l}_{\gamma}} \in \mathfrak{L}_{\mathbf{R}, <\alpha}$. Given $s \in \mathbb{T}_{\mathbf{R}, <\alpha}$, we can write $s$ uniquely as $s = s(\mathfrak{d}_s) \mathfrak{d}_s (1 + \varepsilon_s)$ where $\varepsilon_s \prec 1$. The logarithm on $\mathbb{T}_{\mathbf{R}, <\alpha}$ is then given at $s$ by

$$\log(s) := \log(\mathfrak{d}_s) + \log_{\mathbf{R}}(s(\mathfrak{d}_s)) + \sum_{k > 0} \frac{(-1)^{k+1}}{k} \varepsilon_s^k.$$

See [30, Section 2.3.1] or [6, Section 3.1.3] for more details.

In the case when $\mathbf{R} = \mathbb{R}$ and $\log_{\mathbf{R}}$ is the natural logarithm, we write $\mathfrak{L}_{<\alpha} := \mathfrak{L}_{\mathbb{R}, <\alpha}$ and $\mathbb{L}_{<\alpha} = \mathbb{R}[\![\mathfrak{L}_{<\alpha}]\!]$. The field $\mathbb{L}$ of logarithmic hyperseries is the union

$$\mathbb{L} = \bigcup_{\alpha \in \mathbf{On}} \mathbb{L}_{<\alpha}.$$

This is also a transserial field, and we have $\mathbb{L} = \mathbb{R}[\![\mathfrak{L}]\!]$ where $\mathfrak{L} = \bigcup_{\alpha \in \mathbf{On}} \mathfrak{L}_{<\alpha}$.

**Example 2.3.** The field $\tilde{\mathbb{L}}$ of finitely nested hyperseries can be obtained [9, Section 8] by closing $\mathbb{L}$ under exponential and hyperexponential functions. We write $\tilde{\mathfrak{L}}$ for the monomial group of $\tilde{\mathbb{L}} = \mathbb{R}[\![\tilde{\mathfrak{L}}]\!]$. If $\eta \in \mathbf{On}$, then $\mathbb{L}_{<\omega^{\eta}} = \mathbb{R}[\![\mathfrak{L}_{<\omega^{\nu}}]\!]$ similarly embeds into a minimal closure $\widetilde{\mathbb{L}_{<\omega^{\eta}}} = \mathbb{R}[\![\widetilde{\mathfrak{L}_{<\omega^{\eta}}}]\!]$, and we have $\tilde{\mathbb{L}} = \bigcup_{\eta \in \mathbf{On}} \widetilde{\mathbb{L}_{<\omega^{\eta}}}$.

**Example 2.4.** Écalle's field $\mathbb{T}_{\mathrm{LE}}$ of logarithmic-exponential transseries [19] is *not* a transserial field. Indeed, it is not a field of well-based series, but a directed union $\mathbb{T}_{\mathrm{LE}} = \bigcup_{\mathfrak{G} \in \mathcal{D}} \mathbb{R}[\![\mathfrak{G}]\!]$ of such fields which is a proper subset of the field $\mathbb{R}[\![\bigcup_{\mathfrak{G} \in \mathcal{D}} \mathfrak{G}]\!]$. Although $\mathbb{T}_{\mathrm{LE}}$ is closed under log and its functional inverse exp, none of the fields $\mathbb{R}[\![\mathfrak{G}]\!], \mathfrak{G} \in \mathcal{D}$ are closed under log or exp. Furthermore $\mathbb{T}_{\mathrm{LE}}$ is contained in the class-sized field $\mathbb{R}\langle\!\langle \omega \rangle\!\rangle$ of $\omega$-series, itself a transserial field closed under exp and log.

For the sequel of Section 2, we fix a transserial field $\mathbb{T} = \mathbf{R}[\![\mathfrak{N}]\!]$ with logarithm log.

## 2.2 Exponentiation

The function $\log : \mathbb{T}^{>} \longrightarrow \mathbb{T}$ is strictly increasing, so it has a partially defined left inverse exp, defined by $\exp(\log s) = s$ for all $s \in \mathbb{T}^{>}$. The partial function log is called the *exponential function*. We say that $\mathbb{T}$ is *exponentially closed* if exp is totally defined on $\mathbb{T}$. It can be shown that $\mathbb{T}$ is exponentially closed if and only if

$$\mathbb{T}_{\succ} = \log(\mathfrak{N}). \tag{2.1}$$

If $\mathbb{T}$ is not exponentially closed, then it can be embedded in a canonical way into a minimal exponentially closed transserial field called its exponential closure, see [30, Section 2.3] or [9, Section 8.1]. In particular, the field $\mathbb{R}\langle\!\langle \omega \rangle\!\rangle$ of $\omega$-series is (canonically isomorphic to) the exponential closure of the field $\mathbb{T}_{\mathbb{R}, \log}$ described in Example 2.2.

## 2.3 Powers

Assume that $(\mathbb{T}, \log)$ is exponentially closed, so $\log \colon \mathbb{T}^{>} \longrightarrow \mathbb{T}$ is an isomorphism of ordered groups with inverse exp. As $\mathbb{T}$ is an ordered field extension of $\mathbf{R}$, the ordered group $(\mathbb{T}, +, 0, <)$ is an ordered vector space over $\mathbf{R}$ for scalar multiplication. Therefore $(\mathbb{T}^{>}, \cdot, 1, <)$ is an ordered vector space over $\mathbf{R}$ for the following $\mathbf{R}$-powering operation

$$\begin{aligned} \mathbf{R} \times \mathbb{T}^{>} &\longrightarrow \mathbb{T}^{>} \\ (r, s) &\longmapsto s^r := \exp(r \log(s)). \end{aligned}$$

## 2.4 Steepness

The following ordering on $\mathbb{T}^{\times}$ plays a crucial role in our computations:

**Definition 2.5.** *Let $s, t \in \mathbb{T}^{\times}$. We say that $s$ is **strictly flatter** than $t$ or that $t$ is **strictly steeper** than $s$, and we write.*

$$s \prec\!\!\prec t, \quad \text{if} \quad \log|s| \prec \log|t|.$$

*We say that $s$ is **flatter** than $t$, or that $t$ is **steeper** than $s$, and we write*

$$s \preccurlyeq\!\!\preccurlyeq t, \quad \text{if} \quad \log|s| \preccurlyeq \log|t|.$$

*Lastly, we write*

$$s \asymp\!\!\asymp t \quad \text{if} \quad \log(|s|) \asymp \log(|t|).$$

The relation $\prec\!\!\prec$ is a partial ordering on $\mathbb{T}^{\times}$ with minima $1$ and $-1$.

**Lemma 2.6.** *[30, Proposition 2.2.4] For all $s \in \mathbb{T}^{>, \succ}$, we have $\log s \prec\!\!\prec s$.*

**Lemma 2.7.** *Assume that $(\mathbb{T}, \log)$ is exponentially closed. For all $s, t \in \mathbb{T}^{\times}$, we have*

$$\begin{aligned} s \prec\!\!\prec t &\iff |s|^{\mathbf{R}} < \max(|t|, |t|^{-1}), \\ s \preccurlyeq\!\!\preccurlyeq t &\iff \exists r \in \mathbf{R}^{>}, (\max(|s|, |s|^{-1}))^r < |t|, \\ s \asymp\!\!\asymp t &\iff \exists r_0, r_1 \in \mathbf{R}, |s|^{r_0} < |t| < |s|^{r_1}. \end{aligned}$$

**Proof.** This follows from the relation $\log(|s|^r) = r |s|$ for all $r \in \mathbf{R}$ and the fact that $\log$ is strictly increasing. $\square$

**Corollary 2.8.** *For $s, t \in \mathbb{T}^{>}$, we have*

*a)* $s t \preccurlyeq\!\!\preccurlyeq \underset{\preccurlyeq\!\!\preccurlyeq}{\max}(s, t).$

*b)* $s \asymp t \Longrightarrow s \asymp\!\!\asymp t.$

*c)* $s \prec\!\!\prec t \Longrightarrow s t \asymp\!\!\asymp t.$

**Proof.** The assertions *a)* and *c)* follow from the classical valuation theoretic properties of $\preccurlyeq$. The assertion *b)* is an immediate consequence of Lemma 2.7. $\square$

## 2.5 Transserial derivations

**Definition 2.9.** *A **transserial derivation** on $\mathbb{T}$ is a strongly linear function $\partial \colon \mathbb{T} \longrightarrow \mathbb{T}$ with*

$$\forall s, t \in \mathbb{T}, \partial(s t) = \partial(s) t + s \partial(t) \qquad \text{and} \qquad \forall s \in \mathbb{T}^{\times}, \partial(\log|s|) = \frac{\partial(s)}{s}.$$

*We say that* $(\mathbb{T}, \log, \partial)$ *is a* **transserial differential field**. *An* **identity** *of* $\partial$ *is a series* $x \in \mathbb{T}^{>}$ *with* $1 \notin \operatorname{supp} x$ *and* $\partial(x) = 1$.

**Example 2.10.** The transserial field **No** of surreal numbers, equipped with the Berarducci-Mantova derivation [11], is a transserial differential field with unique identity $\omega$. So are all its differential subfields of the form $\mathbb{R}[\![\mathfrak{M}]\!]$ for a subgroup $\mathfrak{M} \subseteq \mathbf{Mo}$ containing $\omega$. In particular, the class $\mathbb{R}\langle\!\langle \omega \rangle\!\rangle$ of $\omega$-series is also a transserial differential field.

**Example 2.11.** The transserial field $\widetilde{\mathbb{L}}$, together with its standard hyperserial derivation [5, Theorem 5.6], is a transserial differential field with identity $\ell_0$. So are its differential subfields $\widetilde{\mathbb{L}_{<\alpha}}$ for all infinite additively indecomposable ordinals $\alpha$.

**Remark 2.12.** If $\operatorname{Ker}(\partial) = \mathbf{R}$, then any two $s, t \in \mathbb{T}$ with $\partial(s) = \partial(t) = 1$ differ by a constant, so there is a unique identity of $\partial$.

If $\partial$ is a transserial derivation, then we usually write $s' := \partial(s)$ for all $s \in \mathbb{T}$, and

$$t^\dagger := \frac{\partial(t)}{t}$$

for all $t \in \mathbb{T}^\times$. The series $t^\dagger$ is called the *logarithmic derivative* of $t$.

We say that $(\mathbb{T}, \partial)$ is an *H-field* (see [1, 2]) if $\operatorname{Ker}(\partial) = \mathbf{R}$ and $f' > 0$ for all $f \in \mathbb{T}^{>, \succ}$. It is known that **No** (whence also its differential subfield $\mathbb{R}\langle\!\langle \omega \rangle\!\rangle$) is an H-field, as a consequence of [11, Theorem A]. Likewise, the standard derivation on $\widetilde{\mathbb{L}}$ is an H-field [5, Corollary 5.21].

**Lemma 2.13.** *If* $(\mathbb{T}, \partial)$ *is an H-field, then for all* $s, t \in \mathbb{T}$ *we have*

$$t \not\asymp 1 \wedge s \prec t \Longrightarrow s' \prec t'.$$

**Lemma 2.14.** *Assume that* $(\mathbb{T}, \partial)$ *is an H-field. Given and* $\mathfrak{m}, \mathfrak{n} \in \mathfrak{N}$ *with* $\mathfrak{n} \neq 1$, *we have* $\mathfrak{m} \preceq\!\!\preceq \mathfrak{n}$ *if and only if* $\mathfrak{m}^\dagger \preccurlyeq \mathfrak{n}^\dagger$.

**Proof.** Since $\operatorname{supp} \log \mathfrak{n} \succ 1$, we have $\log \mathfrak{n} \not\asymp 1$. The field $(\mathbb{T}, \partial)$ is an H-field, so

$$\mathfrak{m} \preceq\!\!\preceq \mathfrak{n} \Longleftrightarrow \log \mathfrak{m} \preccurlyeq \log \mathfrak{n} \Longleftrightarrow (\log \mathfrak{m})' \preccurlyeq (\log \mathfrak{n})' \Longleftrightarrow \mathfrak{m}^\dagger \preccurlyeq \mathfrak{n}^\dagger. \qquad \square$$

**Definition 2.15.** *Assume that* $(\mathbb{T}, \log, \partial)$ *is an H-field and let* $x \in \mathbb{T}$ *be the identity of* $\partial$. *We say that* $\partial$ *is* **flat** *if for all* $s \in \mathbb{T}$, *all* $\mathfrak{m} \in \mathfrak{N}$ *and all* $\mathfrak{n} \in \operatorname{supp} \mathfrak{m} \circ s$, *there are* $\mathfrak{f}, \mathfrak{w} \in \mathfrak{N}$ *with* $\mathfrak{f} \prec\!\!\prec \mathfrak{m}$ *and* $\mathfrak{w} \preceq\!\!\preceq x$ *and* $\mathfrak{n} = \mathfrak{m} \mathfrak{f} \mathfrak{w}$.

The condition of flatness is satisfied for fields constructed with a variable $x$ and logarithms (or hyperlogarithms) thereof [17, p 14, § 1], and it is preserved under sensible extensions such as exponential or hyperexponential extensions [5, Theorem 5.6].

**Example 2.16.** It is known that the Berarducci-Mantova derivation on **No** is flat [10, Proposition 6.11]. The standard derivation on $\widetilde{\mathbb{L}}$ is flat as a consequence of [5, Theorem 5.6].

**Lemma 2.17.** *Assume that* $\partial$ *is flat. Then for all* $\mathfrak{m} \in \mathfrak{N}$ *with* $\mathfrak{m} \succ x$ *or* ($\mathfrak{m} \prec 1$ *and* $\mathfrak{m} \succcurlyeq x$), *we have* $\operatorname{supp} \mathfrak{m}' \asymp \mathfrak{m}$.

**Proof.** Let $\mathfrak{n} \in \operatorname{supp} \mathfrak{m}'$. By flatness of $\partial$, there are an $\mathfrak{f} \in \mathfrak{N}$ with $\mathfrak{f} \prec\!\!\prec \mathfrak{m}$ and an $\mathfrak{w} \preceq\!\!\preceq x$ with $\mathfrak{w} \preccurlyeq 1$ and $\mathfrak{n} = \mathfrak{m} \mathfrak{f} \mathfrak{w}$. Suppose that $\mathfrak{m} \succ x$. Since $\mathfrak{m} \succ x$, we have $\mathfrak{n} \asymp \mathfrak{m}$ by Corollary 2.8(c).

Suppose that $\mathfrak{m} \succcurlyeq x$ and $\mathfrak{m} \prec 1$. By flatness of $\partial$, for all $\mathfrak{n} \in \operatorname{supp} \mathfrak{m}'$, there is are an $\mathfrak{f} \in \mathfrak{N}$ with $\mathfrak{f} \prec\!\!\prec \mathfrak{m}$ and a $\mathfrak{w} \preceq\!\!\preceq x$ with $\mathfrak{w} \preccurlyeq 1$ and $\mathfrak{n} = \mathfrak{m} \mathfrak{f} \mathfrak{w}$. We have $\mathfrak{n} \preceq\!\!\preceq x$ by Corollary 2.8(a). Now $\log(\mathfrak{m}) + \log(\mathfrak{f}) \sim \log(\mathfrak{m}) < 0$, so $\log(\mathfrak{m}\mathfrak{f}) < 0$, so $\mathfrak{m}\mathfrak{f} \prec 1$. If $\mathfrak{w} \succcurlyeq x$, then it follows since $\mathfrak{w} \preccurlyeq 1$ that $\mathfrak{n} \succcurlyeq x$, hence $\mathfrak{n} \asymp x \asymp \mathfrak{m}$. Otherwise $\mathfrak{w} \prec\!\!\prec x$, so $\mathfrak{n} \asymp \mathfrak{m}$ by Corollary 2.8(c). $\qquad \square$

## 2.6 Transserial composition laws

**Definition 2.18.** *A **transserial composition law** on $\mathbb{T}$ is a function $\circ : \mathbb{T} \times \mathbb{T}^{>,\succ} \longrightarrow \mathbb{T}$ satisfying the following:*

**Compatibility.** *For all $s \in \mathbb{T}^{>,\succ}$, the function $\circ_s : \mathbb{T} \longrightarrow \mathbb{T}; f \mapsto f \circ s$ is a strongly linear morphism of ordered rings which commutes with* log.

**Associativity.** *For all $s, t \in \mathbb{T}^{>,\succ}$ and $f \in \mathbb{T}$, we have $f \circ (s \circ t) = (f \circ s) \circ t$.*

*We say that an element $x \in \mathbb{T}^{>,\succ}$ is an **identity** of $\circ$ if for all $s \in \mathbb{T}^{>,\succ}$ and $f \in \mathbb{T}$, we have $x \circ s = s$ and $f \circ x = f$.*

**Remark 2.19.** If $\circ$ has an identity $x$ then it is unique. Indeed, if $y \in \mathbb{T}^{>,\succ}$ is an identity, then $x = y \circ x = x \circ x$, so by injectivity of $\circ_x$, we have $x = y$.

**Example 2.20.** Given $\eta \in \mathbf{On}$, the standard composition laws on $\mathbb{L}_{<\omega^\eta}$ and $\widetilde{\mathbb{L}_{<\omega^\eta}}$ are transserial composition laws with identity $\ell_0$; see [17, Theorem 1.3] and [5, Result C].

**Definition 2.21.** *A transserial composition law $\circ$ on $\mathbb{T}$ is said **flat** if for all $s \in \mathbb{T}$, all $\mathfrak{m} \in \mathfrak{N}$ and all $\mathfrak{n} \in \operatorname{supp} \mathfrak{m} \circ s$, there are $\mathfrak{f}, \mathfrak{w} \in \mathfrak{N}$ with $\mathfrak{f} \prec\!\!\prec \mathfrak{m}$ and $\mathfrak{w} \preccurlyeq 1$ such that there is an $\mathfrak{u} \in \operatorname{supp} s$ with $\mathfrak{u} \succcurlyeq \mathfrak{w}$, and that $\mathfrak{n} = \mathfrak{d}(\mathfrak{m} \circ s) \mathfrak{f} \mathfrak{w}$.*

Similarly to flat derivations, the condition of flat composition laws is sensible for transserial fields. The standard composition law on $\tilde{\mathbb{L}}$ is flat by [5, Theorems 6.29 and 6.6].

## 2.7 Taylor expansions and chain rule

Consider a transserial derivation $\partial$ and a transserial composition law $\circ$ on $\mathbb{T}$.

**Definition 2.22.** *Let $s \in \mathbb{T}^{>,\succ}$, $f, \delta \in \mathbb{T}$. We say that $f$ **is in Taylor configuration** at $(s, \delta)$ if $\delta \prec s$ and $(\mathfrak{m}' \circ s) \delta \prec \mathfrak{m} \circ s$ for all $\mathfrak{m} \in \operatorname{supp} f$. We say that $f$ **has a Taylor expansion** at $(s, \delta)$ if it is in Taylor configuration and the following holds:*

**TE.** *The family $((\mathfrak{m}^{(k)} \circ s) \delta^k)_{\mathfrak{m} \in \operatorname{supp} f \wedge k \in \mathbb{N}}$ is summable, with*

$$f \circ (s + \delta) = \sum_{k \in \mathbb{N}} \frac{f^{(k)} \circ s}{k!} \delta^k = \sum_{\mathfrak{m} \in \operatorname{supp} f \wedge k \in \mathbb{N}} \frac{\mathfrak{m}^{(k)} \circ s}{k!} \delta^k,$$

*and we have $f \circ s \succ (f' \circ s) \delta \succ \cdots \succ (f^{(k)} \circ s) \delta^k \succ \cdots$.*

We say that *$f$ has Taylor expansions* (with respect to $(\partial, \circ)$) if it has a Taylor expansion at $(s, \delta)$ whenever it is in Taylor configuration at $(s, \delta)$. We say that *$(\mathbb{T}, \partial, \circ)$ has Taylor expansions* if each $f \in \mathbb{T}$ has Taylor expansions.

**Remark 2.23.** By [10, Theorem 3], the structure $(\mathbb{R}\langle\!\langle \omega \rangle\!\rangle, \partial, \circ)$ has Taylor expansions. By [5, Result D], the structure $(\tilde{\mathbb{L}}, \partial, \circ)$ has Taylor expansions.

**Definition 2.24.** *We say that $(\partial, \circ)$ satisfies the **chain rule** if for all $f \in \mathbb{T}$ and $s \in \mathbb{T}^{>,\succ}$, we have $(f \circ s)' = s' f' \circ s$.*

The derivations and composition laws on $\mathbb{R}\langle\!\langle \omega \rangle\!\rangle$ and $\tilde{\mathbb{L}}$ satisfy the chain rule [5, Corollary 6.27].

# 3 Compositional inverses

Let $\mathbf{R}$ be an ordered field. Let $\mathfrak{N}$ be a non-trivial, linearly ordered Abelian group and write $\mathbb{T} = \mathbf{R}[\![\mathfrak{N}]\!]$. Suppose that $(\mathbb{T}, \log, \partial, x)$ is an exponentially closed differential transserial field with identity $x$, that $\partial$ is flat and that $(\mathbb{T}, \partial)$ is an H-field. Let $\circ : \mathbb{T} \times \mathbb{T}^{>,\succ} \longrightarrow \mathbb{T}$ be a flat transserial composition law with identity $x$, and assume that $(\circ, \partial)$ satisfies the chain rule and has Taylor expansions. Let $\mathfrak{M} \subseteq \mathfrak{N}$ be a subgroup containing $x$, and let $\$ = \mathbf{R}[\![\mathfrak{M}]\!]$ be a differential subfield of $\mathbb{T}$ such that

$$\mathfrak{m} \circ (x + \delta) \in \$ \quad \text{whenever } \mathfrak{m} \in \mathfrak{M} \text{ and } \delta \in \$ \text{ with } \delta \prec x. \tag{3.1}$$

Note that $(\$, \partial)$ is an H-field. We will prove the following:

**Theorem 3.1.** *The class $\mathcal{G} := \{x + \delta \in \$ : \delta \prec x\}$ is a group under composition.*

**Lemma 3.2.** *Let $\mathcal{G}_0 \subseteq \$^{>,\succ}$ be such that $(\mathcal{G}_0, \circ, x)$ is a group. Then $(\mathcal{G}_0, \circ, x, <)$ is a linearly right-ordered group.*

**Proof.** The ordering $<$ no $\$$ is linear. For $f, g, h \in \mathcal{G}_0$ with $f < g$, since the function $\circ_h$ is a morphism of ordered rings defined on a linearly ordered field, it is strictly increasing, so $f \circ h < g \circ h$. □

**Lemma 3.3.** *For $f, g \in \$^{>,\succ}$, we have $f \circ g = x \iff g \circ f = x$.*

**Proof.** Assume that $f \circ g = x$. Then $x \circ g = g = g \circ x = g \circ (f \circ g) = (g \circ f) \circ g$. By injectivity of $\circ_g$, we get $x = g \circ f$. The equivalence follows by symmetry. □

We will write $f^{\text{inv}}$ for the inverse of a series $f$ in $(\$^{>,\succ}, \circ)$ when it exists.

**Remark 3.4.** Vincenzo L. Mantova also proposed an independent proof of the fact that $\mathbb{R}\langle\!\langle \omega \rangle\!\rangle^{>,\succ}$ is a group under composition. Our proof of Theorem 3.1 is more general and constructive, but it is based on similar ideas as Mantova's.

## 3.1 Steepness classes

Our method requires us to work with representatives for the steepness equivalence relation $\asymp$. For convenience, we set

$$\mathfrak{P} := \{e^{\mathfrak{n}} : \mathfrak{n} \in \mathfrak{N} \wedge \log x \prec \mathfrak{n}\} \subseteq \mathfrak{N}.$$

Note that for all $\mathfrak{m} \in \mathfrak{M}$ with $\mathfrak{m} \succ x$, we have $\mathfrak{m} \asymp e^{\partial \log \mathfrak{m}}$ where $e^{\partial \log \mathfrak{m}} \in \mathfrak{P}$. Furthermore, for $\mathfrak{n}_0, \mathfrak{n}_1 \in \mathfrak{N}$ with $\log x \prec \mathfrak{n}_0 \prec \mathfrak{n}_1$, we have $x \ll e^{\mathfrak{n}_0} \ll e^{\mathfrak{n}_1}$. So $\mathfrak{P}$ is a class of unique representatives in $\mathfrak{N}$ in each equivalence class in $\{\mathfrak{m} \in \mathfrak{N} : \mathfrak{m} \succ x\}$ for $\asymp$.

**Lemma 3.5.** *Given $\mathfrak{p} \in \mathfrak{P} \cup \{x\}$ and $r \in \mathbf{R}^>$, we have $\mathfrak{p}^\dagger \prec \mathfrak{p}^r$.*

**Proof.** This is immediate if $\mathfrak{p} = x$. Assume that $\mathfrak{p} \succ x$. We have $\log \mathfrak{p} \ll \mathfrak{p}$ by Lemma 2.6. Furthermore, the flatness of $\partial$ entails that $\mathfrak{p}^\dagger = (\log \mathfrak{p})' \preceq \log \mathfrak{p}$, so $\mathfrak{p}^\dagger \ll \mathfrak{p}$. We have $\mathfrak{p}^r \asymp \mathfrak{p}$ by Lemma 2.7 and $\mathfrak{p} \succ 1$, so it follows again by Lemma 2.7 that $\mathfrak{p}^\dagger \prec \mathfrak{p}^r$. □

Let $\mathbf{C} \subseteq e^{\mathfrak{N}}$ be a convex subclass. We write $\mathfrak{M}_{\mathbf{C}} := \{\mathfrak{m} \in \mathfrak{M} : \exists \mathfrak{p}_0, \mathfrak{p}_1 \in \mathbf{C}, \mathfrak{p}_0 \preceq \mathfrak{m} \preceq \mathfrak{p}_1\}$. If $\mathfrak{p}_0 \in \mathfrak{P} \cup \{x\}$, then we write $\mathfrak{M}_{\preceq \mathfrak{p}_0}$ (resp. $\mathfrak{M}_{\prec \mathfrak{p}_0}$, resp. $\mathfrak{M}_{\succ \mathfrak{p}_0}$, resp. $\mathfrak{M}_{\succeq \mathfrak{p}_0}$) for the class $\mathfrak{M}_{\mathbf{C}}$ where $\mathbf{C} = \{\mathfrak{p} \in e^{\mathfrak{N}} : \mathfrak{p} \preceq \mathfrak{p}_0\}$ (resp. $\mathbf{C} = \{\mathfrak{p} \in e^{\mathfrak{N}} : \mathfrak{p} \succeq \mathfrak{p}_0\}$, resp. $\mathbf{C} = \{\mathfrak{p} \in e^{\mathfrak{N}} : \mathfrak{p} \prec \mathfrak{p}_0\}$, resp. $\mathbf{C} = \{\mathfrak{p} \in e^{\mathfrak{N}} : \mathfrak{p} \succ \mathfrak{p}_0\}$). Note that $\mathfrak{M}_{\mathbf{C}}$ is a subgroup of $\mathfrak{M}$, whence $\$_{\mathbf{C}} := \mathbf{R}[\![\mathfrak{M}_{\mathbf{C}}]\!]$ is a subfield of $\$$, whenever $\mathbf{C}$ is an initial subclass of $e^{\mathfrak{N}}$.

**Lemma 3.6.** *Assume that* $\mathbf{C}$ *is a an initial subclass of* $e^{\mathfrak{N}}$ *containing* $x$. *Then the field* $\$_{\mathbf{C}}$ *is a differential subfield of* $(\$, \partial)$.

**Proof.** By strong linearity of $\partial$, it suffices to show that $\partial(\mathfrak{M}_{\mathbf{C}}) \subseteq \$_{\mathbf{C}}$. This follows from Lemma 2.17. □

Likewise, we define $\mathcal{G}_{\mathbf{C}}$ to be the class of series $x + \delta \in \mathcal{G}$ with $\operatorname{supp} \varepsilon \subseteq \mathfrak{M}_{\mathbf{C}}$, and use similar notations $\mathcal{G}_{\preccurlyeq \mathfrak{p}_0}$ for $\mathbf{C} = \{\mathfrak{p} \in \mathfrak{P} : \mathfrak{p} \preccurlyeq \mathfrak{p}_0\}$, and so on. For $\mathfrak{p} \in \mathfrak{P}$, we further write

$$\mathcal{G}_{\mathfrak{p}} := \{x + \varepsilon : \varepsilon \prec 1 \wedge \operatorname{supp} \varepsilon \succcurlyeq \mathfrak{p}\}, \qquad (3.2)$$

and we set

$$\mathcal{G}_x := \mathcal{G}_{\preccurlyeq x} = \{x + \delta : \delta \prec x \wedge \operatorname{supp} \delta \preccurlyeq x\}.$$

## 3.2 Decompositions along steepness classes

Given a series $\varepsilon \in \$^{\prec}$ with $\varepsilon \succ\!\!\!\succ x$ and an $\mathfrak{p} \in \mathfrak{P}$, we write $\varepsilon_{\mathfrak{p}} := \sum_{\mathfrak{m} \succcurlyeq \mathfrak{p}} \varepsilon(\mathfrak{m}) \, \mathfrak{m}$. So each $x + \varepsilon_{\mathfrak{p}}$ for $\mathfrak{p} \in \mathfrak{P}$ lies in $\mathcal{G}_{\mathfrak{p}}$, and $x + \varepsilon = x + \sum_{\mathfrak{p} \in \mathfrak{P}} \varepsilon_{\mathfrak{p}}$. We define

$$\operatorname{supp}_{\mathfrak{P}}(\varepsilon) := \{\mathfrak{p} \in \mathfrak{P} : \varepsilon_{\mathfrak{p}} \neq 0\}.$$

Note that $(\operatorname{supp}_{\mathfrak{P}}(\varepsilon), <)$ embeds into $(\operatorname{supp} \varepsilon, \succ)$ by sending $\mathfrak{p}$ to the least monomial $\mathfrak{m} \in \operatorname{supp} \varepsilon_{\mathfrak{p}}$ with $\mathfrak{p} \asymp \mathfrak{m}$. So the order type $\Lambda(x + \varepsilon)$ of $(\operatorname{supp}_{\mathfrak{P}} \varepsilon, <)$ is an ordinal.

We compute $(x + \varepsilon)^{\operatorname{inv}}$ by writing $x + \varepsilon$ as a well-ordered, right-trailing composition

$$x + \varepsilon = (x + \varepsilon_{\mathfrak{p}_0}) \circ \cdots \circ (x + \varepsilon_{\mathfrak{p}_\gamma}) \circ \cdots, \gamma < \Lambda(x + \varepsilon).$$

The inverse is then obtained as a well-ordered, left-trailing composition

$$\cdots \circ (x + \varepsilon_{\mathfrak{p}_\gamma})^{\operatorname{inv}} \circ \cdots \circ (x + \varepsilon_{\mathfrak{p}_0})^{\operatorname{inv}} = (x + \varepsilon)^{\operatorname{inv}}.$$

We believe there is a well-behaved formal notion of infinite composition on certain groups of formal series, that should be amenable to those simple tricks. We studied an axiomatic framework for such transfinite non-commutative operations [7] and showed that they are well-behaved in the cases of formal series in powers of a variable $x$. Here, we will focus on a minimal treatment of these compositions, and we will only define left-trailing compositions.

## 3.3 Taylor expansions for steepness classes

**Lemma 3.7.** *For* $\mathfrak{p} \in \mathfrak{P} \cup \{x\}$ *and* $x + \varepsilon \in \mathcal{G}_{\succcurlyeq \mathfrak{p}}$, *each* $s \in \$_{\preccurlyeq \mathfrak{p}}$ *has a Taylor expansion at* $(x, \varepsilon)$.

**Proof.** Since $\operatorname{supp} \varepsilon \succcurlyeq \mathfrak{p}$, we have $\varepsilon \preccurlyeq \mathfrak{p}^{-r}$ for some $r \in \mathbf{R}^{>}$. So for all $\mathfrak{m} \in \mathfrak{M}_{\preccurlyeq \mathfrak{p}}$, we have $(\mathfrak{m}^\dagger \circ x)\varepsilon \asymp \mathfrak{m}^\dagger \varepsilon \preccurlyeq \mathfrak{p}^\dagger \varepsilon \preccurlyeq \mathfrak{p}^\dagger \mathfrak{p}^{-r}$. We conclude with Lemma 3.5 that $(\mathfrak{m}^\dagger \circ x)\varepsilon \prec 1$. □

**Corollary 3.8.** *For* $\mathfrak{p} \in \mathfrak{P} \cup \{x\}$ *and* $x + \varepsilon \in \mathcal{G}_{\mathfrak{p}}$, *the function* $\circ_{x + \varepsilon} : \$ \longrightarrow \$$ *restricts to a function* $\$_{\preccurlyeq \mathfrak{p}} \longrightarrow \$_{\preccurlyeq \mathfrak{p}}$ *and* $\circ_{x + \varepsilon} - \operatorname{Id}_{\$}$ *is contracting on this space.*

**Proof.** In view of Lemma 3.7 and **TE**, the function $\circ_{x + \varepsilon} - \operatorname{Id}_{\$}$ is contracting on $\$_{\preccurlyeq \mathfrak{p}}$. For $\mathfrak{m} \in \mathfrak{M}_{\preccurlyeq \mathfrak{p}}$, we have $\mathfrak{m} \circ (x + \varepsilon) = \sum_{k \in \mathbf{N}} \frac{\mathfrak{m}^{(k)}}{k!} \varepsilon^k$. Since $\$$ is closed under $\partial$ and by (3.1), this series lies in $\$$. Now $\operatorname{supp} \mathfrak{m}^{(k)} \preccurlyeq \mathfrak{p}$ by Lemma 3.6 whereas $\operatorname{supp} \varepsilon^k \subseteq (\operatorname{supp} \varepsilon)^k \preccurlyeq \mathfrak{p}$ for all $k \in \mathbf{N}$. So $\mathfrak{m} \circ (x + \varepsilon) \in \$_{\preccurlyeq \mathfrak{p}}$. It follows that $\circ_{x + \varepsilon} - \operatorname{Id}_{\$}$ stabilizes $\$_{\preccurlyeq \mathfrak{p}}$. □

**Proposition 3.9.** *Let* $\mathfrak{p} \in \mathfrak{P} \cup \{x\}$. *Any* $x + \varepsilon \in \mathcal{G}_{\mathfrak{p}}$ *has an inverse in* $\mathcal{G}_{\mathfrak{p}}$.



**Proof.** We claim that the inverse of $x+\varepsilon$ is the series

$$(x+\varepsilon)^{\mathrm{inv}} := (\sum_{i\in\mathbb{N}} (-1)^i \phi^{[i]})(x) \in \mathbb{S}_{\preccurlyeq\mathfrak{p}}$$

where $\phi$ is the restriction of $\circ_{x+\varepsilon} - \mathrm{Id}$ to $\mathbb{S}_{\preccurlyeq\mathfrak{p}}$. Corollary 1.5 gives $(x+\varepsilon)^{\mathrm{inv}} \circ (x+\varepsilon) = x$, so

$$((x+\varepsilon) \circ (x+\varepsilon)^{\mathrm{inv}}) \circ (x+\varepsilon) = (x+\varepsilon) \circ ((x+\varepsilon)^{\mathrm{inv}} \circ (x+\varepsilon)) = x+\varepsilon = x \circ (x+\varepsilon).$$

We deduce since $\circ_{x+\varepsilon}$ is injective that $(x+\varepsilon) \circ (x+\varepsilon)^{\mathrm{inv}} = x$, so $(x+\varepsilon)^{\mathrm{inv}}$ is the inverse of $x+\varepsilon$, as claimed. Since $\phi$ is contracting, we have

$$(x+\varepsilon)^{\mathrm{inv}} - x = \sum_{k>0} (-1)^k \phi^{[k]}(x) \preccurlyeq \phi(x) \prec x.$$

Note that $\phi(x) = \sum_{k>0} \frac{x^{(k)} \circ (x+\varepsilon)}{k!} \varepsilon^k = \varepsilon$. We have $\mathrm{supp}\,((x+\varepsilon)^{\mathrm{inv}} - x) \preccurlyeq \mathfrak{p}$ since $(x+\varepsilon)^{\mathrm{inv}} - x \in \mathbb{S}_{\preccurlyeq\mathfrak{p}}$. If $\mathfrak{p} = x$, then this concludes the proof since $\mathcal{G}_x = \mathcal{G}_{\preccurlyeq x}$. If $\mathfrak{p} \succ x$, then since $\mathrm{supp}\,((x+\varepsilon)^{\mathrm{inv}} - x) \preccurlyeq \varepsilon$, we have $\mathrm{supp}\,((x+\varepsilon)^{\mathrm{inv}} - x) \succcurlyeq \mathfrak{p}$, so $(x+\varepsilon)^{\mathrm{inv}} \in \mathcal{G}_{\mathfrak{p}}$. □

**Lemma 3.10.** *For* $\mathfrak{p} \in \mathfrak{P}$, $f, g \in \mathcal{G}_{\preccurlyeq\mathfrak{p}}$ *with* $f \circ g = x$, *and* $\iota \in \mathbb{S}_{\succ\mathfrak{p}}$. *Write* $h := f \circ (g+\iota)$, *so* $h \in \mathbb{S}$ *by* (3.1). *Then*

a) $\mathrm{supp}_{\mathfrak{P}}(h-x) = \mathrm{supp}_{\mathfrak{P}} \iota$.

b) *For all* $\mathfrak{t} \in \mathrm{supp}_{\mathfrak{P}} \iota$, *we have* $\max_{\prec}\{\tau \in \mathrm{term}\,\iota : \tau \asymp \mathfrak{t}\} = \max_{\prec}\{\tau \in \mathrm{term}\,h : \tau \asymp \mathfrak{t}\}$.

**Proof.** Write $\varepsilon := f - x$ and $\delta := g - x$. Let $\mathfrak{m} \in \mathrm{supp}\,\varepsilon$. The chain rule gives $(\mathfrak{m} \circ g)^{\dagger} = g' \mathfrak{m}^{\dagger} \circ g$. Since $(\mathbb{S}, \partial)$ is an H-field, we have $g' \sim 1$, so $(\mathfrak{m} \circ g)^{\dagger} \asymp \mathfrak{m}^{\dagger} \circ g$. We have $\mathfrak{m}^{\dagger} \prec \mathfrak{m}$ by Lemma 3.5, so $\mathfrak{m}^{\dagger} \in \mathbb{S}_{\preccurlyeq\mathfrak{p}}$. Recall by Lemma 3.7 that $\mathfrak{m}^{\dagger}$ has a Taylor expansion at $(x, \delta)$. In particular, we have $\mathfrak{m}^{\dagger} \circ g \sim \mathfrak{m}^{\dagger}$. We deduce by Lemma 3.7 that $\mathfrak{m}$ is in Taylor configuration at $(g, \iota)$. By **TE**, we obtain

$$h - x = f \circ (g+\iota) - x = f \circ (g+\iota) - f \circ g = \sum_{k>0} \frac{f^{(k)} \circ g}{k!} \iota^k. \qquad (3.3)$$

Let $k > 0$. The chain rule gives $f' \circ g = \frac{1}{g'}$. We deduce by induction that there is a rational fraction $F_k(X_1, \ldots, X_k) \in \mathbb{Q}(X_1, \ldots, X_k)$ with $f^{(k)} \circ g = F_k(g', \ldots, g^{(k)})$. By Corollary 2.8(a), we deduce that $\mathrm{supp}\,f^{(k)} \circ g \preccurlyeq \mathfrak{p}$. Therefore Corollary 2.8(c) gives $\frac{\mathrm{supp}\,f^{(k)} \circ g\, \iota^k}{\iota^k} \prec \iota^k$, so

$$\mathrm{supp}\,(f' \circ g\, \iota) \succ \mathrm{supp}\,(f^{(2)} \circ g\, \iota^2) \succ \cdots.$$

We get

$$\begin{aligned}
\mathrm{supp}_{\mathfrak{P}}(f \circ (g+\iota) - f \circ g) &= \bigcup_{k>0} \mathrm{supp}_{\mathfrak{P}}(f^{(k)} \circ g\, \iota^k) \\
&= \bigcup_{k>0} \mathrm{supp}_{\mathfrak{P}} \iota^k \qquad \text{(by Corollary 2.8(c))} \\
&= \mathrm{supp}_{\mathfrak{P}} \iota.
\end{aligned}$$

This proves a. Let $\mathfrak{t} \in \mathrm{supp}_{\mathfrak{P}} \iota$. Write $\tau_0 := \max_{\prec}\{\tau \in \mathrm{term}\,\iota : \tau \asymp \mathfrak{t}\}$ and $\tau_1 := \max_{\prec}\{\tau \in \mathrm{term}\,h : \tau \asymp \mathfrak{t}\}$. Note that for all infinitesimals $\mathfrak{u}, \mathfrak{v} \in \mathfrak{M}$ with $\mathfrak{u} \prec \tau_0 \preccurlyeq \mathfrak{v}$, we have $\mathfrak{u} \succ \mathfrak{v}$ so $\mathfrak{u}\tau_0 \succ \mathfrak{v}\tau_0$. In view of (3.3), we have $\tau_1 = \frac{\tau_{f' \circ g}}{1} \tau_0 = \tau_{g'} \tau_0$. But $g = x + \varepsilon$ for an $\varepsilon \prec x$, so since $(\mathbb{T}, \partial)$ is an H-field, we have $\varepsilon' \prec 1$ and $\tau_{g'} = 1$, hence the result. □

**Proposition 3.11.** *If* $\mathcal{G}_{\succ x}$ *is a group, then* $\mathcal{G}$ *is a group.*



**Proof.** Assume that $\mathcal{G}_{\succ x}$ is a group. Note that its identity is $x$. Let $f = x + \delta \in \mathcal{G}$. We can write $\delta$ uniquely as $\delta = \delta_{\preccurlyeq x} + \delta_{\succ x}$ where $\delta_{\preccurlyeq x} \in \mathbb{S}_{\preccurlyeq x}$ and $\delta_{\succ x} \in \mathbb{S}_{\succ x}$. The series $x + \delta_{\preccurlyeq x}$ has an inverse $g$ in $\mathcal{G}_x$ by Lemma 3.9, and we have $g \circ f \in \mathcal{G}_{\succ x}$ by Lemma 3.10(a). Let $h = x + \varepsilon$ be the inverse of $g \circ f$ in $\mathcal{G}_{\succ x}$, so $h \circ g$ is the inverse of $f$ in $\mathbb{S}^{>,\succ}$. Now $h \circ g = g + \varepsilon \circ g$ where $\varepsilon \prec 1$ and $g - x \prec x$, so $h \circ g \in \mathcal{G}$. □

## 3.4 Well-ordered left-trailing compositions

Let $\alpha$ be an ordinal, let $(\mathfrak{p}_\beta)_{\beta < \alpha}$ be a strictly increasing sequence in $\mathfrak{P}$, and let $\mathcal{E} := (x + \varepsilon_{\mathfrak{p}_\beta})_{\beta < \alpha}$ be a sequence with $x + \varepsilon_{\mathfrak{p}_\beta} \in \mathcal{G}_{\mathfrak{p}_\beta}$ for each $\beta < \alpha$. We define a composition

$$"\cdots \circ (x + \varepsilon_{\mathfrak{p}_\beta}) \circ \cdots \circ (x + \varepsilon_{\mathfrak{p}_0})" = \mathrm{LC}(\mathcal{E})$$

by induction on $\beta \leqslant \alpha$. Set $\mathrm{LC}(\mathcal{E} \upharpoonright 0) := x$. If $\beta = \gamma + 1$ is a successor, then set

$$(x + \varepsilon_{\mathfrak{p}_\gamma}) \circ \mathrm{LC}(\mathcal{E} \upharpoonright \gamma) =: \mathrm{LC}(\mathcal{E} \upharpoonright \beta). \tag{3.4}$$

Finally, if $\beta$ is a non-zero limit, set

$$\mathrm{LC}(\mathcal{E} \upharpoonright \beta) := x + \sum_{\gamma < \beta} \mathrm{LC}(\mathcal{E} \upharpoonright \gamma + 1) - \mathrm{LC}(\mathcal{E} \upharpoonright \gamma). \tag{3.5}$$

We first justify that (3.5) is warranted, i.e. that the family $(\mathrm{LC}(\mathcal{E} \upharpoonright \gamma + 1) - \mathrm{LC}(\mathcal{E} \upharpoonright \gamma))_{\gamma < \beta}$ is summable. Note that when that family is defined and summable, we have

$$\mathrm{LC}(\mathcal{E} \upharpoonright \eta) = \mathrm{LC}(\mathcal{E} \upharpoonright \alpha) + \sum_{\alpha \leqslant \gamma < \eta} \mathrm{LC}(\mathcal{E} \upharpoonright \gamma + 1) - \mathrm{LC}(\mathcal{E} \upharpoonright \gamma) \tag{3.6}$$

for all $\alpha < \eta \leqslant \beta$. In view of (3.4), for each $\alpha < \eta \leqslant \beta$, this can be written as

$$\mathrm{LC}(\mathcal{E} \upharpoonright \eta) = \mathrm{LC}(\mathcal{E} \upharpoonright \alpha) + \sum_{\alpha \leqslant \gamma < \eta} \varepsilon_\gamma \circ \mathrm{LC}(\mathcal{E} \upharpoonright \gamma). \tag{3.7}$$

**Lemma 3.12.** *Suppose that the family* $(\mathrm{LC}(\mathcal{E} \upharpoonright \gamma + 1) - \mathrm{LC}(\mathcal{E} \upharpoonright \gamma))_{\gamma < \beta}$ *is well-defined and summable. For all $\eta < \beta$, writing $\mathcal{E}^{\uparrow \eta}$ for the function defined on the unique ordinal $\mu$ with $\eta + \mu = \beta$ by $\forall \iota < \mu, \mathcal{E}^{\uparrow \eta}(\iota) = \mathcal{E}(\eta + \iota)$, we have*

$$\mathrm{LC}(\mathcal{E} \upharpoonright \beta) = \mathrm{LC}(\mathcal{E}^{\uparrow \eta}) \circ \mathrm{LC}(\mathcal{E} \upharpoonright \eta).$$

**Proof.** We prove this by induction on $\beta$. This is vacuously true if $\beta = 0$ and it is immediate by induction and in view of (3.4) if $\beta$ is a successor. Assume that $\beta$ is a non-zero limit and that the result holds for strictly smaller ordinals. We have

$$\begin{aligned}
\mathrm{LC}(\mathcal{E}^{\uparrow \eta}) \circ \mathrm{LC}(\mathcal{E} \upharpoonright \eta) &= \left( x + \sum_{\iota < \mu} \mathrm{LC}(\mathcal{E}^{\uparrow \eta} \upharpoonright \iota + 1) - \mathrm{LC}(\mathcal{E}^{\uparrow \eta} \upharpoonright \iota) \right) \circ \mathrm{LC}(\mathcal{E} \upharpoonright \eta) \\
&= \mathrm{LC}(\mathcal{E} \upharpoonright \eta) + \sum_{\iota < \mu} \mathrm{LC}(\mathcal{E}^{\uparrow \eta} \upharpoonright \iota + 1) \circ \mathrm{LC}(\mathcal{E} \upharpoonright \eta) - \mathrm{LC}(\mathcal{E}^{\uparrow \eta} \upharpoonright \iota) \circ \mathrm{LC}(\mathcal{E} \upharpoonright \eta) \\
&= \mathrm{LC}(\mathcal{E} \upharpoonright \eta) + \sum_{\iota < \mu} \mathrm{LC}(\mathcal{E} \upharpoonright \eta + \iota + 1) - \mathrm{LC}(\mathcal{E} \upharpoonright \eta + \iota) \quad \text{(ind. hyp.)} \\
&= \mathrm{LC}(\mathcal{E} \upharpoonright \eta) + \sum_{\eta \leqslant \gamma < \beta} \mathrm{LC}(\mathcal{E} \upharpoonright \gamma) - \mathrm{LC}(\mathcal{E} \upharpoonright \gamma) \\
&= \mathrm{LC}(\mathcal{E} \upharpoonright \beta).
\end{aligned}$$

The result follows by induction. □



**Lemma 3.13.** *Suppose that* $(\operatorname{LC}(\mathcal{E}\restriction \gamma+1)-\operatorname{LC}(\mathcal{E}\restriction \gamma))_{\gamma<\beta}$ *is well-defined and summable. Then for all* $\mathfrak{n}\in\operatorname{supp}\operatorname{LC}(\mathcal{E}\restriction \beta)$*, there are a* $\gamma<\beta$ *and a* $\mathfrak{q}\in\operatorname{supp}\varepsilon_\gamma$ *with* $\mathfrak{n}\preccurlyeq\mathfrak{q}\circ\operatorname{LC}(\mathcal{E}\restriction \gamma)$.

**Proof.** We prove the result by induction on $\beta$. It holds trivially with $(k,\gamma_k,\mathfrak{n}_k)=(1,0,x)$ if $\beta=0$. Let $\beta>0$ such that the result holds for all $\gamma\in\beta$. In view of (3.7), there is a $\gamma<\beta$ with $\mathfrak{n}\in\operatorname{supp}\varepsilon_\gamma\circ\operatorname{LC}(\mathcal{E}\restriction \gamma)$. So $\mathfrak{n}\in\operatorname{supp}\mathfrak{m}_0\circ\operatorname{LC}(\mathcal{E}\restriction \gamma)$ for a certain $\mathfrak{m}_0\in\operatorname{supp}\varepsilon_\gamma$. Let $\mathfrak{n}_0$ be the dominant monomial of $\mathfrak{m}_0\circ\operatorname{LC}(\mathcal{E}\restriction \gamma)$. Since $\circ$ is flat, there are an $\mathfrak{w}\in\operatorname{supp}\operatorname{LC}(\mathcal{E}\restriction \gamma)$ and an $\mathfrak{f}\prec\mathfrak{n}_0$ with $\frac{\mathfrak{n}}{\mathfrak{n}_0}=\mathfrak{f}\mathfrak{w}$. If $\mathfrak{n}_0\succcurlyeq\mathfrak{w}$, then $\mathfrak{n}\preccurlyeq\mathfrak{n}_0\asymp\mathfrak{m}_0\circ\operatorname{LC}(\mathcal{E}\restriction \gamma)$ by Corollary 2.8(a), so we are done. Otherwise, we have $\mathfrak{n}\preccurlyeq\mathfrak{w}$, so the induction hypothesis yields the result. $\square$

**Lemma 3.14.** *Let* $\beta\leqslant\alpha$ *be an ordinal such that each* $\operatorname{LC}(\mathcal{E}\restriction \gamma)$ *for* $\gamma<\beta$ *is defined. Then for* $\eta<\gamma<\beta$*, and* $\mathfrak{q}\in\operatorname{supp}\varepsilon_\eta$*, we have* $\operatorname{supp}\mathfrak{q}\circ\operatorname{LC}(\mathcal{E}\restriction \eta)\prec\!\prec\varepsilon_\gamma\circ\operatorname{LC}(\mathcal{E}\restriction \gamma)$.

**Proof.** We proceed by induction. If $\beta\leqslant 1$, then this is vacuously true. Let $1<\beta$ such that the result holds for all $\alpha<\beta$, for all families $\mathcal{E}$ of domain $\eta$. Let $\eta<\gamma<\beta$ and let $\mathfrak{q}\in\operatorname{supp}\varepsilon_\eta$. We first show that

$$\mathfrak{q}\circ\operatorname{LC}(\mathcal{E}\restriction \eta)\prec\!\prec\varepsilon_\gamma\circ(\operatorname{LC}(\mathcal{E}\restriction \gamma)). \tag{3.8}$$

We have $\operatorname{LC}(\mathcal{E}\restriction \gamma)=\operatorname{LC}(\mathcal{E}^{\uparrow\eta}\restriction \iota)\circ\operatorname{LC}(\mathcal{E}\restriction \eta)$, so it suffices to show that

$$\mathfrak{q}\prec\!\prec\varepsilon_\gamma\circ\operatorname{LC}(\mathcal{E}^{\uparrow\eta}\restriction \iota). \tag{3.9}$$

The induction hypothesis for the family $\mathcal{E}^{\uparrow\eta}$ gives $\operatorname{supp}\operatorname{LC}(\mathcal{E}^{\uparrow\eta}\restriction \iota)-x\succcurlyeq\varepsilon_\eta$, so Lemma 3.7 yields $\mathfrak{q}\sim\mathfrak{q}\circ\operatorname{LC}(\mathcal{E}^{\uparrow\eta}\restriction \iota)$, whence (3.9) reduces to the valid inequality $\mathfrak{q}\prec\!\prec\varepsilon_\gamma$.

Now let $\mathfrak{n}\in\operatorname{supp}\mathfrak{q}\circ\operatorname{LC}(\mathcal{E}\restriction \eta)$. By flatness of $\circ$, we have $\frac{\mathfrak{n}}{\mathfrak{q}\circ\operatorname{LC}(\mathcal{E}\restriction \eta)}\preccurlyeq\mathfrak{w}$ for some $\mathfrak{w}\in\operatorname{supp}\operatorname{LC}(\mathcal{E}\restriction \eta)$. By Lemma 3.13, we have $\mathfrak{w}\preccurlyeq\mathfrak{u}\circ\operatorname{LC}(\mathcal{E}\restriction \alpha)$ for some $\alpha<\eta$ and $\mathfrak{u}\in\operatorname{supp}\varepsilon_\alpha$. The induction hypothesis gives $\mathfrak{u}\circ\operatorname{LC}(\mathcal{E}\restriction \alpha)\prec\!\prec\mathfrak{q}\circ\operatorname{LC}(\mathcal{E}\restriction \eta)$, so $\mathfrak{n}\asymp\mathfrak{q}\circ\operatorname{LC}(\mathcal{E}\restriction \eta)$, whence $\mathfrak{n}\prec\!\prec\varepsilon_\gamma\circ(\operatorname{LC}(\mathcal{E}\restriction \gamma))$ by (3.8). This concludes the proof. $\square$

**Corollary 3.15.** *We have* $\operatorname{LC}(\mathcal{E}\restriction \gamma)\triangleleft\operatorname{LC}(\mathcal{E}\restriction \gamma+1)$ *for all* $\gamma<\alpha$.

We conclude that $(\operatorname{LC}(\mathcal{E}\restriction \gamma+1)-\operatorname{LC}(\mathcal{E}\restriction \gamma))_{\gamma<\alpha}$ is well-defined and summable, so $\operatorname{LC}(\mathcal{E})$ is well-defined and satisfies the above properties. By (3.1) the series $\operatorname{LC}(\mathcal{E})$ lies in $\$$.

## 3.5 Inverses as transfinite compositions

In this subsection, we prove Theorem 3.1 by computing the inverses of all elements of $\mathcal{G}_{\succ x}$. Let $f=x+\varepsilon\in\mathcal{G}_{\succ x}$ and let $(\mathfrak{p}_\gamma)_{\gamma<\Lambda(f)}$ be the unique strictly $\prec\!\prec$-increasing parametrisation of $\operatorname{supp}_{\mathfrak{P}}\varepsilon$. For each $\beta\leqslant\Lambda(f)$, we write $f_\beta:=x+\sum_{\gamma<\beta}\varepsilon_{\mathfrak{p}_\gamma}$, and $\delta_\beta:=f-f_\beta$. We will define a family $\mathcal{E}=(x+\varepsilon_\gamma)_{\gamma<\Lambda(f)}$ by induction on $\beta\leqslant\Lambda(f)$.

For simplicity, given $\gamma\leqslant\beta\leqslant\Lambda(f)$ and $\sigma\leqslant\beta$ with $\gamma+\sigma=\beta$, we will write $g_\gamma:=\operatorname{LC}(\mathcal{E}\restriction \gamma)$ and $\rho_{\gamma,\beta}:=\operatorname{LC}(\mathcal{E}^{\uparrow\gamma}\restriction \sigma)$. So $g_\beta=\rho_{\gamma,\beta}\circ g_\gamma$ for all $\gamma\leqslant\beta$. We will assume inductively that for all $\beta\leqslant\Lambda(f)$, we have

$$\forall\gamma<\beta, x+\varepsilon_\gamma \ \in \ \mathcal{G}_{\mathfrak{p}_\gamma}, \tag{3.10}$$
$$\forall\gamma\leqslant\beta, g_\gamma \ = \ f_\gamma^{\operatorname{inv}}, \tag{3.11}$$
$$\forall\gamma\leqslant\beta, \operatorname{supp}_{\mathfrak{P}}(g_\gamma\circ f-x) \ = \ \operatorname{supp}_{\mathfrak{P}}f\cap\{\mathfrak{p}\in\mathfrak{P}:\forall\eta<\gamma,\mathfrak{p}\succ\mathfrak{p}_\eta\}. \tag{3.12}$$

Moreover, we make the additional inductive assumption

**I$_\beta$.** Let $\mathfrak{p}\in\operatorname{supp}_{\mathfrak{P}}\varepsilon$ with $\mathfrak{p}\succ\mathfrak{p}_\gamma$ for all $\gamma<\beta$. Write

$$\tau_\gamma(\mathfrak{p}):=\max_{\prec}\{\tau\in\operatorname{term}g_\gamma\circ f:\tau\asymp\mathfrak{p}\}$$

for each $\gamma < \beta$. Then $\tau_{\gamma+1}(\mathfrak{p}) = \tau_\gamma(\mathfrak{p})$ for all $\gamma < \beta$.

If (3.12) holds at $\beta = \Lambda(f)$, then we conclude that $\mathrm{supp}_{\mathfrak{P}}(g_\beta \circ f - x) = \varnothing$, whence $g_\beta \circ f = x$. We prove the result by double induction. Let $\eta$ be an ordinal such that the conditions hold for all series $h \in \mathcal{G}_{\succ x}$ with $\Lambda(h) < \eta$. Consider an $f$ as above with $\Lambda(f) = \eta$, let $\beta \leqslant \Lambda(f)$ such that the conditions on $\mathcal{E}$ are satisfied for all $\gamma < \beta$.

An immediate induction shows that our definition of $g_\beta$ only depends on $f_\beta$. Hence if $f = f_\beta$, then we obtain $g_\beta = f_\beta^{\mathrm{inv}}$ immediately. So (3.11) holds. In order to derive the other statements, we distinguish between the successor and limit cases.

**Successor case.** Suppose that $\beta = \gamma + 1$ is a successor. If $\beta = \Lambda(f)$, then we are done. Otherwise, (3.12) at $\gamma$ implies that $g_\gamma \circ f$ has the form $x + \iota + \delta$ where $x + \iota \in \mathcal{G}_{\mathfrak{p}_\gamma} \setminus \{x\}$ and $\mathrm{supp}_{\mathfrak{P}} \delta = \{\mathfrak{p} \in \mathrm{supp}_{\mathfrak{P}} f : \mathfrak{p} \succ \mathfrak{p}_\gamma\}$. We then set $x + \varepsilon_\gamma := (x + \iota)^{\mathrm{inv}} \in \mathcal{G}_{\mathfrak{p}_\gamma}$, whence (3.10) holds. So

$$g_\beta \circ f = (x + \iota)^{\mathrm{inv}} \circ g_\gamma \circ f = (x + \iota)^{\mathrm{inv}} \circ (x + \iota + \delta).$$

We have $\mathrm{supp}\,(g_\beta \circ f - x) = \mathrm{supp}_{\mathfrak{P}} \delta$ by Lemma 3.10(a), so (3.12) holds at $\beta$.

Finally, let us prove $\mathbf{I}_\beta$. Let $\mathfrak{p} \in \mathrm{supp}_{\mathfrak{P}} \varepsilon$ with $\mathfrak{p} \succ \mathfrak{p}_\eta$ for all $\gamma < \beta$, and let $\eta < \beta$. If $\eta + 1 < \beta$, then we have $\tau_{\eta+1}(\mathfrak{p}) = \tau_\eta(\mathfrak{p})$ by the induction hypothesis. Suppose now that $\eta + 1 = \beta$, so $\eta = \gamma$. Lemma 3.10(b) implies that $\tau_\beta(\mathfrak{p}) = \tau_\gamma(\mathfrak{p})$. So $\mathbf{I}_\beta$ holds.

**Limit case.** Suppose now that $\beta$ is a non-zero limit. Note that (3.10) and $\mathbf{I}_\beta$ immediately follow from the induction hypothesis. Let $\gamma < \beta$ and let $\sigma$ with $\gamma + \sigma = \beta$. Set $h_\gamma := g_\gamma \circ f$. Recall that $g_\gamma = f_{\varphi(\gamma)}^{\mathrm{inv}}$ by (3.11), so $f = f_{\varphi(\gamma)} \circ h_\gamma$.

The relation $g_\beta = \rho_{\gamma,\beta} \circ g_\gamma$ gives $g_\beta \circ f = \rho_{\gamma,\beta} \circ h_\gamma$. By (3.7), we can also write $\rho_{\gamma,\beta} = x + \sum_{\eta < \sigma} \Delta_\eta$, where $\Delta_\eta := \varepsilon_{\mathfrak{p}_{\gamma+\eta}} \circ \mathrm{LC}(\mathcal{E}^{\uparrow \gamma} \upharpoonright \gamma + \eta)$. Thus

$$g_\beta \circ f = h_\gamma + \sum_{\eta < \sigma} \Delta_\eta \circ h_\gamma$$

Fix a $\mathfrak{p} \in \mathrm{supp}_{\mathfrak{P}}(g_\beta \circ f)$. Then $\mathfrak{p} \in \mathrm{supp}_{\mathfrak{P}} h_\gamma$ or $\mathfrak{p} \in \mathrm{supp}_{\mathfrak{P}}(\Delta_\eta \circ h_\gamma)$ for some $\eta < \sigma$.

Assume for contradiction that $\mathfrak{p} \notin \bigcup_{\alpha < \beta} \mathrm{supp}_{\mathfrak{P}} h_\alpha$. By (3.12), there is a $\gamma_0 < \beta$ with $\mathfrak{p} \notin \mathrm{supp}_{\mathfrak{P}} h_\gamma$ whenever $\gamma \geqslant \gamma_0$. Writing $\gamma_0 + \sigma_0 = \beta$, we have

$$\mathfrak{p} \in \bigcap_{\eta < \sigma_0} \mathrm{supp}_{\mathfrak{P}}(\Delta_\eta \circ h_{\gamma_0}). \qquad (3.13)$$

Let $\eta < \sigma_0$ and let $\mathfrak{m} \in \mathrm{supp}\,\Delta_\eta \circ h_{\gamma_0}$ be given, as a consequence of (3.13), with $\mathfrak{p} \asymp \mathfrak{m}$. By flatness of $\circ$, there are an $\mathfrak{n} \in \mathrm{supp}\,\Delta_\eta$, an $\mathfrak{f} \prec \mathfrak{n} \circ h_{\gamma_0}$ and a $\mathfrak{w} \in \mathrm{supp}\,h_{\gamma_0} \cup \{1\}$ such that $\frac{\mathfrak{m}}{(\mathfrak{n} \circ h_{\gamma_0}) \mathfrak{f}} \asymp \mathfrak{w}$. By Lemma 3.14, we have $\mathrm{supp}\,\Delta_\eta \prec \mathrm{supp}\,\Delta_\rho$ for all $\rho > \eta$. Choosing $\eta$ (and $\gamma_0$ if necessary) large enough, we may assume that $\mathfrak{n} \circ h_{\gamma_0} \not\asymp \mathfrak{w}$. If $\mathfrak{n} \circ h_{\gamma_0} \succ \mathfrak{w}$, then $\mathfrak{p} \asymp \mathfrak{n} \circ h_{\gamma_0}$ by Corollary 2.8(c). So again choosing $\eta$ large enough, we may assume that $\mathfrak{n} \circ h_{\gamma_0} \prec \mathfrak{w}$, so $\mathfrak{m} \asymp \mathfrak{w}$ by Corollary 2.8(c). But then $\mathfrak{p} \in \mathrm{supp}_{\mathfrak{P}} h_{\gamma_0}$: a contradiction. Therefore we must have

$$\mathrm{supp}_{\mathfrak{P}}(g_\beta \circ f) \subseteq \bigcap_{\alpha < \beta} \mathrm{supp}_{\mathfrak{P}} h_\alpha \subseteq \mathrm{supp}_{\mathfrak{P}} f \cap \{\mathfrak{t} \in \mathfrak{P} : \forall \gamma < \beta, \mathfrak{t} \succ \mathfrak{p}_\gamma\}.$$

Conversely, consider a $\mathfrak{p} \in \mathrm{supp}_{\mathfrak{P}} f$ with $\mathfrak{p} \succ \mathfrak{p}_\gamma$ for all $\gamma < \beta$. In view of $\mathbf{I}_\beta$ we have $\mathfrak{d}_\tau \notin \mathrm{supp}\,h_{\eta+1} - h_\eta$ for all $\eta < \beta$ and $\mathfrak{d}_\tau \in \mathrm{supp}\,h_0$, so $\mathfrak{d}_\tau \in \mathrm{supp}\,(h_0 + (\sum_{\eta < \beta} h_{\eta+1} - h_\eta)) = \mathrm{supp}\,g_\beta \circ f$, whence $\mathfrak{p} \in \mathrm{supp}_{\mathfrak{P}} g_\beta \circ f$. This shows (3.12).

Together with Lemma 3.3 and Proposition 3.11, this establishes Theorem 3.1.

## 3.6 Applications

Consider the class $\tilde{\mathbb{L}} = \mathbb{R}[\![\tilde{\mathfrak{L}}]\!]$ of finitely nested hyperseries. By Theorem 3.1, we have:

**Theorem 3.16.** *Let $\mathfrak{M} \subseteq \tilde{\mathfrak{L}}$ be a subgroup containing $\ell_0$ and assume that $\mathbb{S} = \mathbb{R}[\![\mathfrak{M}]\!]$ is closed under derivation and composition. Then $(\{\ell_0 + \delta : \delta \in \mathbb{S} \wedge \delta \prec \ell_0\}, \circ, \ell_0)$ is a group.*

**Lemma 3.17.** *Let $f \in \tilde{\mathbb{L}}^{>,\succ}$ and let $\lambda$ be a limit ordinal with $f \in \widetilde{\mathbb{L}_{<\omega^\lambda}}$. There is an infinitesimal $\varepsilon \in \widetilde{\mathbb{L}_{<\omega^\lambda}} \circ \mathrm{e}_{\omega^\lambda}^{\ell_0}$ with $\ell_{\omega^\lambda} \circ f \circ \mathrm{e}_{\omega^\lambda}^{\ell_0} = \ell_0 + \varepsilon$.*

**Proof.** By [5, Lemma 6.24 and eq. (2.2)] there are a $\mu < \lambda$ such that for all $\nu \geqslant \mu$, there is a $\delta \in \widetilde{\mathbb{L}_{<\omega^\lambda}}$ with $\delta \prec \ell_{\omega^\mu}$ and $\ell_{\omega^\mu} \circ f = \ell_{\omega^\mu} + \delta$. Let $\ell_{\omega^\lambda}^{\uparrow \omega^\mu} \in \mathbb{L}$ be the unique series [17, Section 5] with $\ell_{\omega^\lambda}^{\uparrow \omega^\mu} \circ \ell_{\omega^\mu} = \ell_{\omega^\lambda}$, and set $\varepsilon_0 := \sum_{k>0} \frac{\left(\ell_{\omega^\lambda}^{\uparrow \omega^\mu}\right)^{(k)} \circ \ell_{\omega^\mu}}{k!} \delta^k$ By **TE**, we have

$$\ell_{\omega^\lambda} \circ f - \ell_{\omega^\mu} = \ell_{\omega^\lambda}^{\uparrow \omega^\mu} \circ (\ell_{\omega^\mu} + \delta) - \ell_{\omega^\lambda} = \varepsilon_0,$$

Each $(\ell_{\omega^\lambda}^{\uparrow \omega^\mu})^{(k)}$ for $k > 0$ lies in $\mathbb{L}_{<\omega^\lambda}$, so $\varepsilon_0 \in \widetilde{\mathbb{L}_{<\omega^\lambda}}$. So $\varepsilon := \varepsilon_0 \circ \mathrm{e}_{\omega^\lambda}^{\ell_0}$ satisfies the conditions. □

**Corollary 3.18.** *The class $(\tilde{\mathbb{L}}^{>,\succ}, \circ, \ell_0)$ is a group.*

We next turn to the class $\mathbb{R}\langle\!\langle \omega \rangle\!\rangle$ of $\omega$-series.

**Lemma 3.19.** *Let $a \in \mathbb{R}\langle\!\langle \omega \rangle\!\rangle^{>,\succ}$. There are an $(e, n) \in \mathbb{Z} \times \mathbb{N}$ and an $\varepsilon \in \mathbb{R}\langle\!\langle \omega \rangle\!\rangle^{\prec}$ with*

$$\log_n(\omega) \circ (a \circ \log_e(\omega)) \circ \exp_n(\omega) = \omega + \varepsilon.$$

**Proof.** It suffices to show that there is an $(e, n) \in \mathbb{Z} \times \mathbb{N}$ with $\log_n(\omega) \circ a \circ \exp_n(\omega) - \exp_e(\omega) \prec 1$ (then compose on the right by $\log_e(\omega)$). This follows from [11, Corollary 5.9]. □

**Corollary 3.20.** *The class $(\mathbb{R}\langle\!\langle \omega \rangle\!\rangle^{>,\succ}, \circ, \omega)$ is a group.*

**Proof.** The class $(\{\omega + \delta : \delta \in \mathbb{R}\langle\!\langle \omega \rangle\!\rangle \wedge \delta \prec \omega\}, \circ, \omega)$ is a group by Theorem 3.1. Since each $\log_e(\omega), e \in \mathbb{Z}$ is invertible in $(\mathbb{R}\langle\!\langle \omega \rangle\!\rangle^{>,\succ}, \circ)$, the result follows by Lemma 3.19. □

## 4 Conjugacy

In [20], Écalle studies what he calls the natural growth scale. This is a (somewhat informally defined) group $\mathcal{G}$, under composition, of germs at $+\infty$ belonging in quasianalytic classes of real-valued functions. The elements in $\mathcal{G}$ involve transexponential and sublogarithmic functions $\exp_{\omega^k}, \log_{\omega^k}, k \in \mathbb{N}$ which satisfy the same conjugacy equations

$$\begin{aligned}
\exp_{\omega^{k+1}}(r+1) &= \exp_{\omega^k}(\exp_{\omega^{k+1}}(r)) \quad \text{and} \\
\log_{\omega^{k+1}}(\log_{\omega^k}(r)) &= \log_{\omega^{k+1}}(r)
\end{aligned}$$

of $\ell_{\omega^k}$ and $\mathrm{e}_{\omega^k}^{\ell_0}$, for large enough $r \in \mathbb{R}$ (see [30, Appendix A]).

The groups $(\tilde{\mathbb{L}}^{>,\succ}, \circ)$ and $(\mathcal{G}, \circ)$ can be regarded respectively as formal and geometric substantiations of the same idea. Écalle gives formulas for conjugacy relations within $(\mathcal{G}, \circ)$. In order to make sense of those formulas in our formal setting, we rely on Edgar's work [21]. Edgar shows that each transseries $f \in \mathbb{T}_{\mathrm{LE}}$ with $f > x$ of exponentiality 0 is a conjugate of $x + 1$ [21, Theorem 4.4]. Edgar's proofs apply in the case of $\tilde{\mathbb{L}}$, as well as for more general series fields, with a few adjustments that will be made below.

**Setting:** Throughout Section 4, we fix a structure $(\mathbb{T}, \log, \partial, \circ, x)$ where

1. $\mathbb{T} = \mathbf{R}[\![\mathfrak{N}]\!]$ where $\mathfrak{N}$ is a linearly ordered Abelian group.

2. $(\mathbf{R}, +, \cdot, 0, 1, <, \exp)$ is a model of the first-order theory of the real ordered exponential field $\mathbb{R}_{\exp}$.

3. $(\mathbb{T}, \log, \partial)$ is an exponentially closed transserial field.

4. $\partial : \mathbb{T} \longrightarrow \mathbb{T}$ is a *surjective* transserial derivation with identity $x$, which is flat.

5. $(\mathbb{T}, \partial)$ is an H-field.

6. $\circ : \mathbb{T} \times \mathbb{T}^{>,\succ} \longrightarrow \mathbb{T}$ is a transserial composition law with identity $x$, which is flat.

7. $(\mathbb{T}, \partial, \circ)$ has Taylor expansions and satisfies the chain rule.

8. $(\mathbb{T}^{>,\succ}, \circ, x)$ is a group.

**Remark 4.1.** Condition 2 on $\mathbf{R}$ is necessary. Indeed Edgar's method relies on the resolution of differential equations [21, Section 1, (4)] involving coefficient functions $\alpha_{\mathfrak{m}} : \mathbb{R} \longrightarrow \mathbb{R}$ for $\mathfrak{m}$ ranging in certain subgroups of $\mathfrak{M}$. More precisely, writing $\mathcal{C}^{\infty}(\mathbb{R})$ for the algebra of smooth real-valued functions, he makes use of the derivation

$$\begin{aligned} \mathrm{d} : \mathcal{C}^{\infty}(\mathbb{R})[\![\mathfrak{M}]\!] &\longrightarrow \mathcal{C}^{\infty}(\mathbb{R})[\![\mathfrak{M}]\!] \\ \sum_{\mathfrak{m} \in \mathfrak{M}} \alpha_{\mathfrak{m}} \mathfrak{m} &\longmapsto \sum_{\mathfrak{m} \in \mathfrak{M}} \alpha'_{\mathfrak{m}} \mathfrak{m}. \end{aligned}$$

on the algebra $\mathcal{C}^{\infty}(\mathbb{R})[\![\mathfrak{M}]\!]$ of functions $\alpha : \mathfrak{M} \longrightarrow \mathcal{C}^{\infty}(\mathbb{R})$ ; $\mathfrak{m} \mapsto \alpha_{\mathfrak{m}}$ with well-based support. The two relevant properties of d are its surjectivity and its kernel $\mathrm{Ker}(\mathrm{d}) = \mathbb{T}$. These fail if one replaces $\mathcal{C}^{\infty}(\mathbb{R})$ with the algebra of infinitely differentiable functions $\mathbf{R} \longrightarrow \mathbf{R}$ for arbitrary ordered fields $\mathbf{R}$ (see [6, Section 2.1]). Yet, assuming that $\mathbf{R}$ satisfies 2, they hold if $\mathcal{C}^{\infty}(\mathbb{R})$ is replaced by the algebra $\mathbf{R}[x, \exp]$ of functions of the form $\mathbf{R} \longrightarrow \mathbf{R}$; $r \mapsto P(r, \exp(r))$ for $P \in \mathbf{R}[X, Y]$. Such functions are definable in $(\mathbf{R}, +, \cdot, 0, 1, <, \exp)$. Thus the fact that $\mathrm{d} \upharpoonright \mathbf{R}[x, \exp] : \mathbf{R}[x, \exp] \longrightarrow \mathbf{R}[x, \exp]$ is surjective with kernel $\mathbf{R}$ can be translated into a set of first-order statements in the language of ordered exponential rings. As $\mathbf{R}$ and $\mathbb{R}_{\exp}$ are elementarily equivalent, the map d is surjective with kernel $\mathbb{T}$. One easily shows in view of [21, (17)] that the functions $\alpha_{\mathfrak{m}}$ are all in $\mathbf{R}[x, \exp]$. Thus all of Edgar's computations can be carried out identically in our context.

We further assume the following monotonicity condition:

**Monotonicity.** For all $f \in \mathbb{T} \setminus \mathbf{R}$, the function

$$\begin{aligned} \hat{f} : \mathbb{T}^{>,\succ} &\longrightarrow \mathbb{T} \\ s &\longmapsto f \circ s \end{aligned}$$

is strictly increasing if $f' > 0$ and strictly decreasing if $f' < 0$.

**Remark 4.2.** All the conditions holds when $\mathbb{T} = \tilde{\mathbb{L}}$. For 4, this is a consequence of [5, Result B]. For 8, this is Corollary 3.18. For monotonicity, this is [5, Result E]. The other conditions are already checked in Section 2.

**Proposition 4.3.** *The structure* $(\mathbb{T}^{>,\succ}, \circ, x, <)$ *is a (linearly) bi-ordered group.*

**Proof.** We already know by Lemma 3.2 that $(\mathbb{T}^{>,\succ}, \circ, x, <)$ is a (linearly) right-ordered group. Let $f, g, h \in \mathbb{T}^{>,\succ}$ with $f > g$. By monotonicity, the function $\hat{h}$ is strictly increasing, so $h \circ f > h \circ g$, i.e. $(\mathbb{T}^{>,\succ}, \circ, x, <)$ is left-ordered. □

**Lemma 4.4.** *For* $f, s, t \in \mathbb{T}^{>,\succ}$ *with* $s \neq t$, *we have* $f \circ t - f \circ s \sim \tau_f \circ t - \tau_f \circ s$.

**Proof.** The result is trivial if $s=t$, so we may assume without loss of generality that $t>s$. Write $\delta = f - \tau_f$, so $\delta \prec f$. Let $r \in \mathbf{R}$. The series $g := f - r\,\delta$ lies in $\mathbb{T}^{>,\succ}$, so $g \circ t > g \circ s$, i.e. $f \circ t - r\,\delta \circ t > f \circ s - r\,\delta \circ s$, so $f \circ t - f \circ s > r\,(\delta \circ t - \delta \circ s)$. We deduce that $f \circ t - f \circ s \succ \delta \circ t - \delta \circ s$, i.e. $f \circ t - f \circ s \sim \tau_f \circ t - \tau_f \circ s$. $\square$

We fix a subgroup $\mathfrak{M} \subseteq \mathfrak{N}$ containing $x$ such that $\mathbb{S} := \mathbf{R}[\![\mathfrak{M}]\!]$ is a differential subfield of $\mathbb{T}$ with $\partial(\mathbb{S}) = \mathbb{S}$, that $\mathbb{S}$ is closed under $\circ$, and that $(\mathbb{S}^{>,\succ}, \circ, x)$ is a group.

**Definition 4.5.** *A **steep subfield** of $\mathbb{S}$ is a differential subfield $\mathbb{U} := \mathbf{R}[\![\mathfrak{U}]\!]$ of $\mathbb{S}$ where $\mathfrak{U} \subseteq \mathfrak{M}$ is a subgroup containing $x$ and satisfying*

**Steepness.** *For all $\mathfrak{u} \in \mathfrak{U}^{\neq}$, we have $\mathfrak{u} \succeq x$.*

*A **steep conjugate** of $\mathbb{S}$ is a steep subfield $\mathbb{U}$ of $\mathbb{S}$ satisfying*

**Conjugacy.** *For all $f \in \mathbb{S}^{>,\succ}$ with $f > x$, there are a $V \in \mathbb{S}^{>,\succ}$ and an $\varepsilon \in \mathbb{U}^{\prec}$ with $\varepsilon > 0$ and $V \circ f = (x + \varepsilon) \circ V$.*

Our main general theorem is as follows:

**Theorem 4.6.** *Assume that $\mathbb{S}$ has a steep conjugate. Then for all $f \in \mathbb{S}^{>,\succ}$ with $f > x$, there is a $V \in \mathbb{S}^{>,\succ}$ with $V \circ f = V + 1$.*

Theorem 4.6 follows from the following proposition.

**Proposition 4.7.** *Let $\mathbb{U}$ be a steep subfield of $\mathbb{S}$. For all $\varepsilon \in \mathbb{U}^{\prec}$ with $\varepsilon > 0$, there is a $V \in \mathbb{S}^{>,\succ}$ with $V \circ (x + \varepsilon) = V + 1$.*

Suppose for a moment that $\mathbb{S}$ has a steep conjugate and that Theorem 4.6 holds. We see by taking inverses that any two $f, g \in \mathbb{S}^{>,\succ}$ with $f, g < x$ are conjugate. No two series $f, g \in \mathbb{S}^{>,\succ}$ with $f < x < g$ can be conjugate, so there are exactly three conjugacy classes in $(\mathbb{S}^{>,\succ}, \circ)$, including $\{x\}$. In particular, the group $(\mathbb{S}^{>,\succ}, \circ)$ is simple. Furthermore, the positive cone $\mathbb{S}^{>x}$ (hence the ordering on $\mathbb{S}^{>,\succ}$) is first-order definable with parameters in $(\mathbb{S}^{>,\succ}, \circ)$ as the class of series that are conjugates of $x + 1$. Throughout the proof of Proposition 4.7, we fix a steep subfield $\mathbb{U} = \mathbf{R}[\![\mathfrak{U}]\!]$ of $\mathbb{S}$.

## 4.1 Edgar's method

Besides the properties of $\mathbb{T}_{\mathrm{LE}}$ as a differential field (in particular, that it is an H-field), Edgar relies on properties of an integration operator $\int$ on $\mathbb{T}_{\mathrm{LE}}$. So we must introduce it in our setting. We have $\partial(\mathbb{T}) = \mathbb{T}$ by assumption and $\mathrm{Ker}(\partial) = \mathbf{R}$ since $(\mathbb{T}, \partial)$ is an H-field. So for each $f \in \mathbb{T}$, there is a unique $F \in \mathbb{T}$ with $1 \notin \mathrm{supp}\, F$ and $F' = f$. We write $\int f := F$. For $s, t \in \mathbb{T}^{>,\succ}$, we also write

$$\int_s^t f := \left(\int f\right) \circ t - \left(\int f\right) \circ s.$$

**Lemma 4.8.** *For $s, t \in \mathbb{T}^{>,\succ}$, the functions $\int : \mathbb{T} \longrightarrow \mathbb{T}$ and $\int_s^t : \mathbb{T} \longrightarrow \mathbb{T}$ are strongly linear.*

**Proof.** By Lemma 2.13, given $\mathfrak{m} \in \mathfrak{M}$, there is a unique $\mathfrak{n} \in \mathfrak{M}$ and a unique $r \in \mathbf{R}^{\times}$ with $\mathfrak{m} \sim r\,\mathfrak{n}'$. We then write $\mathcal{I}(\mathfrak{m}) := r\,\mathfrak{n}$. Note that $\mathcal{I} : \mathfrak{M} \longrightarrow \mathfrak{M}$ is strictly $\prec$-increasing so it extends uniquely into a strongly linear function $\mathcal{I} : \mathbb{T} \longrightarrow \mathbb{T}$. The strongly linear map

$$\begin{aligned} \phi : \mathbb{T} &\longrightarrow \mathbb{T} \\ f &\longmapsto f - (\partial \circ \mathcal{I})(f) \end{aligned}$$

is contracting, i.e. satisfies $\phi(\mathfrak{m}) \prec \mathfrak{m}$ for all $\mathfrak{m} \in \mathfrak{M}$. By Corollary 1.5, the map $\partial \circ \mathcal{I} = \mathrm{id}_\mathbb{T} - \phi$ has a strongly linear inverse $(\partial \circ \mathcal{I})^{\mathrm{inv}} = \sum_{k \in \mathbb{N}} \Phi^{[k]}$. We have $(\mathcal{I} \circ (\partial \circ \mathcal{I})^{\mathrm{inv}})(1) = \mathcal{I}(0) = 0$. We deduce that $\int = \mathcal{I} \circ (\partial \circ \mathcal{I})^{\mathrm{inv}}$ is strongly linear, whence also $\int_s^t$ is strongly linear. $\square$

We next derive elementary properties of the integral operator.

**Lemma 4.9.** *Let $f, g \in \mathbb{T}$ and $s, t \in \mathbb{T}^{>, \succ}$ with $s \leqslant t$. If $0 \leqslant f \leqslant g$ then $0 \leqslant \int_s^t f \leqslant \int_s^t g$.*

**Proof.** By linearity, it is enough to show that $\int_s^t f \geqslant 0$. Since $(\int f)' = f \geqslant 0$, the function $\mathbb{T}^{>, \succ} \to \mathbb{T}; u \mapsto (\int f) \circ u$ is nondecreasing by monotonicity. Therefore $\int_s^t f = (\int f) \circ t - (\int f) \circ s \geqslant 0$. $\square$

**Corollary 4.10.** *For $f, g \in \mathbb{T}$ and $s, t \in \mathbb{T}^{>, \succ}$ with $s \neq t$, we have $f \prec g \Longrightarrow \int_s^t f \prec \int_s^t g$.*

**Proof.** We may assume that $s < t$. For all $r \in \mathbf{R}^>$, we have $|f| \leqslant r |g|$. Therefore
$$\left| \int_s^t f \right| = \int_s^t |f| \leqslant \int_s^t r |g| = r \int_s^t |g| = r \left| \int_s^t g \right|$$
by Lemmas 4.9 and 4.8. Note that $g \neq 0$, so $\int g \notin \mathbf{R}$. Since $s \neq t$, we have $(\int g) \circ s \neq (\int g) \circ t$ by monotonicity, whence $|\int_s^t g| > 0$. We deduce that $|\int_s^t f| \prec |\int_s^t g|$, whence $\int_s^t f \prec \int_s^t g$. $\square$

**Corollary 4.11.** *For $f, g \in \mathbb{T}^\times$ and $s, t \in \mathbb{T}^{>, \succ}$ with $s \neq t$, we have $f \sim g \Longrightarrow \int_s^t f \sim \int_s^t g$.*

We require some technical lemmas whose proofs in [21] rely on the specific inductive construction of $\mathbb{T}_{\mathrm{LE}}$, and must therefore be proved in a different way. Using the adapted versions of those results, we will prove an adapted version of [21, Theorem 4.4].

**Lemma 4.12.** (adapted from [21, Lemma 3.11(g)]) *If $\mathfrak{b} \in \mathfrak{U}^{\prec}$ and $\mathfrak{n} \in \mathrm{supp}\, (x\,\mathfrak{b})'$, then $\mathfrak{n}^\dagger \preccurlyeq \mathfrak{b}^\dagger$.*

**Proof.** We have $(x\,\mathfrak{b})' = \mathfrak{b} + x\,\mathfrak{b}'$ so we may assume that $\mathfrak{n} \in x \cdot \mathrm{supp}\, \mathfrak{b}'$. If $\mathfrak{b} \succcurlyeq x$, then $\mathfrak{b}^\dagger \asymp x^\dagger$. Let $\mathfrak{m} \in \mathrm{supp}\, \mathfrak{b}'$ with $\mathfrak{n} = x\,\mathfrak{m}$. We have $\mathfrak{m} \succcurlyeq \mathfrak{b}$ by Lemma 2.17, so $\mathfrak{m}^\dagger \asymp \mathfrak{b}^\dagger$, so $\mathfrak{n}^\dagger = x^\dagger + \mathfrak{m}^\dagger \preccurlyeq \mathfrak{b}^\dagger$. Otherwise $\mathfrak{b} \succ x$ by steepness of $\mathbb{U}$. We have $\log((\mathrm{supp}\, \mathfrak{b}') \cdot x) \subseteq (\log(\mathrm{supp}\, \mathfrak{b}')) + \log x$ where $\log(\mathrm{supp}\, \mathfrak{b}') \asymp \log \mathfrak{b} \succ \log x$ by Lemma 2.17, whence $(\log(\mathrm{supp}\, \mathfrak{b}')) \succ \log x$. In particular $\log \mathfrak{n} \sim \log \mathfrak{b}$ so $\mathfrak{n}^\dagger \sim \mathfrak{b}^\dagger$, whence $\mathfrak{n}^\dagger \preccurlyeq \mathfrak{b}^\dagger$. $\square$

**Lemma 4.13.** (adapted from [21, Lemma 3.14(a,b)]) *Let $\mathfrak{m} \in \mathfrak{M}$ and set*
$$\mathfrak{B} := \{\mathfrak{n} \in \mathfrak{M} : \mathfrak{n}^\dagger \prec \mathfrak{m}\} \qquad \text{and} \qquad \overline{\mathfrak{B}} := \{\mathfrak{n} \in \mathfrak{M} : \mathfrak{n}^\dagger \preccurlyeq \mathfrak{m}\}.$$
*Then $\mathfrak{B}$ and $\overline{\mathfrak{B}}$ are subgroups of $\mathfrak{M}$. Moreover, if $\mathfrak{g} \in \mathfrak{U}^{\prec}$, then we have*
$$\mathfrak{g} \in \mathfrak{B} \Longrightarrow \mathrm{supp}\,(x\,\mathfrak{g})' \subseteq \mathfrak{B} \qquad \text{and} \qquad \mathfrak{g} \in \overline{\mathfrak{B}} \Longrightarrow \mathrm{supp}\,(x\,\mathfrak{g})' \subseteq \overline{\mathfrak{B}}.$$

**Proof.** The fact that $\mathfrak{B}$ and $\overline{\mathfrak{B}}$ are subgroups follows from the inequality $(\mathfrak{n}_0\,\mathfrak{n}_1)^\dagger = \mathfrak{n}_0^\dagger + \mathfrak{n}_1^\dagger \preccurlyeq \max_{\preccurlyeq}(\mathfrak{n}_0^\dagger, \mathfrak{n}_1^\dagger)$ for all $\mathfrak{n}_0, \mathfrak{n}_1 \in \mathfrak{M}$. The last two statements follow from Lemma 4.12. $\square$

**Lemma 4.14.** (adapted from [21, Lemma 3.20]) *Let $\mathfrak{B} \subseteq \mathfrak{M}$ be well-based. Let $\mathfrak{g} \in \mathfrak{B}$. There are finitely many pairs $\mathfrak{g}_1, \mathfrak{g}_2 \in \mathfrak{B}$ with $\mathfrak{g} \in \mathrm{supp}\,(x\,\mathfrak{g}_1)'\,\mathfrak{g}_2$.*

**Proof.** The family $(x\,\mathfrak{m})_{\mathfrak{m} \in \mathfrak{B}}$ is summable. Since $\partial$ is strongly linear, it follows that $((x\,\mathfrak{m})')_{\mathfrak{m} \in \mathfrak{B}}$ is summable. So $((x\,\mathfrak{m})'\,\mathfrak{n})_{\mathfrak{m}, \mathfrak{n} \in \mathfrak{B}}$ is summable by [24, Proposition 3.3]. $\square$

**Lemma 4.15.** *(adapted from [21, Lemma 3.21]) Let $\mathfrak{e} \in \mathfrak{M}^{\prec}$, and set*

$$\mathfrak{S} := \{\mathfrak{g} \in \mathfrak{M} : \mathfrak{g} \preccurlyeq \mathfrak{e} \wedge \mathfrak{g}^\dagger \prec (x\,\mathfrak{e})^{-1}\} \quad \text{and} \quad \overline{\mathfrak{S}} := \{\mathfrak{g} \in \mathfrak{M} : \mathfrak{g} \preccurlyeq \mathfrak{e} \wedge \mathfrak{g}^\dagger \preccurlyeq (x\,\mathfrak{e})^{-1}\}.$$

*Then $\mathbf{R}[\![\mathfrak{S}]\!]$ and $\mathbf{R}[\![\overline{\mathfrak{S}}]\!]$ are closed under the operations*

$$\begin{aligned} \mathfrak{g}_0, \mathfrak{g}_1 &\mapsto \mathfrak{g}_0\,\mathfrak{g}_1, \\ \mathfrak{g}_0, \mathfrak{g}_1 &\mapsto (x\,\mathfrak{g}_0)'\,\mathfrak{g}_1, \quad \text{and} \\ \mathfrak{g} &\mapsto x\,\mathfrak{e}\,\mathfrak{g}'. \end{aligned}$$

**Proof.** Apply Edgar's proof of [21, Lemma 3.21], using Lemma 4.13 instead of [21, Lemma 3.14(a,b)]. □

**Lemma 4.16.** *(adapted from [21, Lemma 3.23]) Let $\mathfrak{B} \subseteq \mathfrak{U}$ be non-empty, well-based and infinitesimal. Write $\mathfrak{e} = \max \mathfrak{B}$ and assume that $\mathfrak{g}^\dagger \preccurlyeq (x\,\mathfrak{e})^{-1}$ for all $\mathfrak{g} \in \mathfrak{B}$. Let $\overline{\mathfrak{B}}$ denote the smallest subclass of $\mathfrak{U}$ such that*

  i. $\overline{\mathfrak{B}} \supseteq \mathfrak{B}$,
  ii. *if $\mathfrak{g}_1, \mathfrak{g}_2 \in \overline{\mathfrak{B}}$, then $\mathfrak{g}_1\,\mathfrak{g}_2 \in \overline{\mathfrak{B}}$, and*
  iii. *if $\mathfrak{g}_1, \mathfrak{g}_2 \in \overline{\mathfrak{B}}$, then $\operatorname{supp}(x\,\mathfrak{g}_1)'\,\mathfrak{g}_2 \subseteq \overline{\mathfrak{B}}$.*

*Then $\overline{\mathfrak{B}}$ is well-based.*

**Proof.** We need to prove that the least set $\mathfrak{B}_1$ of monomials with $\mathfrak{B}_1 \supseteq \mathfrak{B} \cup \{\mathfrak{e}^2\}$ with $\forall \mathfrak{g} \in \mathfrak{B}_1, \operatorname{supp}(x\,\mathfrak{e}\,\mathfrak{g}') \subseteq \mathfrak{B}_1$ is well-based. To that end, write

$$\mathfrak{S} := \{\mathfrak{m} \in \mathfrak{U}^{\prec} : \mathfrak{m}^\dagger \prec (x\,\mathfrak{e})^{-1}\} \quad \text{and} \quad \overline{\mathfrak{S}} := \{\mathfrak{m} \in \mathfrak{U}^{\prec} : \mathfrak{m}^\dagger \preccurlyeq (x\,\mathfrak{e})^{-1}\}.$$

Consider the derivation $\partial_\mathfrak{e} := x\,\mathfrak{e}\,\partial$ on $\mathbb{T}$. We claim that

$$\partial_\mathfrak{e}(\mathbf{R}[\![\mathfrak{S}]\!]) \subseteq \mathbf{R}[\![\mathfrak{S}]\!] \quad \text{and} \quad \partial_\mathfrak{e}(\mathbf{R}[\![\overline{\mathfrak{S}}]\!]) \subseteq \mathbf{R}[\![\overline{\mathfrak{S}}]\!].$$

Given $\mathfrak{m} \in \overline{\mathfrak{S}}$, we have $\operatorname{supp} \mathfrak{m}' \precapprox \mathfrak{m}$ by Lemma 2.17. We deduce that $(\operatorname{supp} \mathfrak{m}')^\dagger \asymp \mathfrak{m}^\dagger \preccurlyeq (x\,\mathfrak{e})^{-1}$ with a strict inequality if $\mathfrak{m} \in \mathfrak{S}$. Thus $(\operatorname{supp} \partial_\mathfrak{e}(\mathfrak{m}))^\dagger \subseteq \{x^{-1}, \mathfrak{e}^\dagger\} \cup (\operatorname{supp} \mathfrak{m}')^\dagger \preccurlyeq (x\,\mathfrak{e})^{-1}$, with strict inequality if $\mathfrak{m} \in \mathfrak{S}$. So $\partial_\mathfrak{e}$ restricts to strongly linear maps $\mathbf{R}[\![\mathfrak{S}]\!] \longrightarrow \mathbf{R}[\![\mathfrak{S}]\!]$ and $\mathbf{R}[\![\overline{\mathfrak{S}}]\!] \longrightarrow \mathbf{R}[\![\overline{\mathfrak{S}}]\!]$. Moreover we have $\partial_\mathfrak{e}(\mathfrak{m}) \prec \mathfrak{m}$ for all $\mathfrak{m} \in \mathfrak{S}$. So we have a strictly extensive and Noetherian choice operator on $\mathbf{R}[\![\mathfrak{S}]\!]$, as per [23, Appendix A.4], given by $\vartheta(\mathfrak{m}) := \operatorname{supp} \partial_\mathfrak{e}(\mathfrak{m})$ for all $\mathfrak{m} \in \mathfrak{S}$. Given $\mathfrak{X} \subseteq \mathfrak{U}$, we write $\operatorname{Cl}(\mathfrak{X})$ for the union of classes $\operatorname{Cl}_n(\mathfrak{X}), n \in \mathbb{N}$ where

$$\begin{aligned} \operatorname{Cl}_0(\mathfrak{X}) &:= \mathfrak{X} \quad \text{and} \\ \operatorname{Cl}_{n+1}(\mathfrak{X}) &:= \operatorname{Cl}_n(\mathfrak{X}) \cup \bigcup_{\mathfrak{m} \in \operatorname{Cl}_n(\mathfrak{X})} \operatorname{supp} \partial_\mathfrak{e}(\mathfrak{m}) \quad \text{for all } n \in \mathbb{N}. \end{aligned}$$

As a consequence of [3, Corollary 1.4], for each well-based subset $\mathfrak{W}$ of $\mathfrak{S}$, the set $\operatorname{Cl}(\mathfrak{W})$ is well-based. Let $\mathfrak{C} := (\mathfrak{B} \cup \{\mathfrak{e}^2\}) \cap \mathfrak{S}$ and $\mathfrak{D} := (\mathfrak{B} \cup \{\mathfrak{e}^2\}) \setminus \mathfrak{S} = \{\mathfrak{m} \in \mathfrak{B} \cup \{\mathfrak{e}^2\} : \mathfrak{m}^\dagger \asymp (x\,\mathfrak{e})^{-1}\}$. Since $\partial_\mathfrak{e}(\mathfrak{S}) \subseteq \mathbf{R}[\![\mathfrak{S}]\!]$, writing $\mathfrak{E} = \bigcup_{\mathfrak{m} \in \mathfrak{D}} \operatorname{supp} \partial_\mathfrak{e}(\mathfrak{m}) \setminus \{\mathfrak{m}\}$, we have $\mathfrak{E} \subseteq \mathfrak{S}$ and $\operatorname{Cl}(\mathfrak{D}) = \mathfrak{D} \cup \operatorname{Cl}(\mathfrak{E})$. So $\mathfrak{B}_1 = \operatorname{Cl}(\mathfrak{B}) \subseteq \operatorname{Cl}(\mathfrak{C}) \cup \operatorname{Cl}(\mathfrak{D}) \subseteq \operatorname{Cl}(\mathfrak{C}) \cup \mathfrak{D} \cup \operatorname{Cl}(\mathfrak{E})$ is well-based. Now apply Edgar's proof of [21, Lemma 3.23], using Lemma 4.15 instead of [21, Lemma 3.21(c)]. □

## 4.2 Solving conjugacy equations

**Proposition 4.17.** *(adapted from [21, Theorems 3.8 and 4.1]) Let $f \in \mathbb{U}^{>, \succ}$ of the form $f = x\,(1 + r\,\mathfrak{e} + \delta)$ where $r \in \mathbf{R}^{>}$, $\delta \prec \mathfrak{e} \prec 1$, $(\operatorname{supp} \delta)^\dagger \preccurlyeq (x\,\mathfrak{e})^{-1}$. Then there is a $V \in \mathbb{S}^{>, \succ}$ with $V' \sim (r\,x\,\mathfrak{e})^{-1}$ and $V \circ f = V + 1$.*

**Proof.** Note that $\mathfrak{e}, \delta \in \mathbb{U}$. By the arguments in the proof of [21, Theorems 3.8], using Lemma 4.16 instead of [21, Lemma 3.23] and Lemma 4.14 instead of [21, Lemma 3.20], we obtain a series $\Phi_1(0, x) \in \mathbb{S}$ with $\Phi_1(0, x) \sim r\, x\, \mathfrak{e}$ and

$$\frac{f'}{\Phi_1(0, x) \circ f} = \frac{1}{\Phi_1(0, x)}.$$

Set $V := \int \frac{1}{\Phi_1(0,x)} \in \int \mathbb{S}$. We have $(V' \circ f)\, f' = V'$, whence $(V \circ f)' = V'$ by the chain rule. So $V \circ f = V + r_0$ for a certain $r_0 \in \mathbf{R}$. It is enough in order to conclude to prove that $r_0 = 1$. Set $A = \int \frac{1}{r\, x\, \mathfrak{e}}$ and write $\tau$ for the dominant term of $A$. So

$$V \circ f - V = \int_x^f \frac{1}{\Phi_1(0, x)} \sim \int_x^f \frac{1}{r\, x\, \mathfrak{e}}, \qquad \text{(by Corollary 4.11)}$$
$$\sim A \circ f - A \sim \tau \circ f - \tau. \qquad \text{(by Lemma 4.4)}$$

We have $\frac{1}{r\, x\, \mathfrak{e}} \succ x^{-1} = (\log x)'$ so $A \succ 1$. It follows that $f - x \asymp \mathfrak{e}\, x \prec \frac{\tau}{\tau'}$. Thus by **TE**, we have $\tau \circ f - \tau \sim \tau'(f - x) \sim A'(f - x) \sim \frac{r\, x\, \mathfrak{e}}{r\, x\, \mathfrak{e}} \sim 1$. So $r_0 \sim 1$. But $r_0 \in \mathbf{R}$, so $r_0 = 1$. □

**Proposition 4.18.** *(adapted from [21, Theorem 4.2]) Let $f \in x + \mathbb{U}^\prec$ where $f - x > 0$. There are a $V \in \mathbb{S}^{>, \succ}$ and a $\delta \in \mathbb{S}^\prec$ with $\delta^\dagger \succ 1$ and $V \circ f \circ V^{\mathrm{inv}} = x + 1 + \delta$.*

**Proof.** The proof is the same as in [21, Proposition 4.2], using Proposition 4.17 instead of [21, Theorem 4.1]. □

**Proposition 4.19.** *(adapted from [21, Proposition 4.3]) Let $f \in \mathbb{S}^{>, \succ}$ of the form $f = x + 1 + \delta$ where $\delta \prec 1$ and $\delta^\dagger \succ 1$. There is a $V \in \mathbb{S}^{>, \succ}$ with*

$$V \circ f = V + 1. \tag{4.1}$$

**Proof.** Write $\Gamma$ for the function $\mathbb{S} \longrightarrow \mathbb{S}$ defined by $\Gamma(g) := (x + g) \circ f - x - 1 = g \circ f + \delta$. It suffices to show that $\Gamma$ has a fixed point $g$; then $V := x + g$ satisfies (4.1).

To that end, we will show that there is a subclass $\mathfrak{S}$ of $\mathfrak{M}$ such that $\Gamma(\mathbf{R}[\![\mathfrak{S}]\!]) \subseteq \mathbf{R}[\![\mathfrak{S}]\!]$ and that $\mathfrak{m} \circ f \prec \mathfrak{m}$ for all $\mathfrak{m} \in \mathfrak{S}$. Consider the class $\mathfrak{S}$ of monomials $\mathfrak{m}$ with $\mathfrak{m} \prec 1$ and $\mathfrak{m}^\dagger \succ 1$. We have $\delta \in \mathbf{R}[\![\mathfrak{S}]\!]$ since $\mathfrak{d}_\delta^\dagger \succ 1$ and the logarithmic derivative is $\prec$-nonincreasing on $\mathfrak{N}^\prec$. Fix an $\mathfrak{m} \in \mathfrak{S}$. We have

$$\log \frac{\mathfrak{m} \circ f}{\mathfrak{m}} = \log(\mathfrak{m} \circ f) - \log \mathfrak{m}$$
$$= \int_x^f \mathfrak{m}^\dagger$$
$$\succ \int_x^f 1. \qquad \text{(by Corollary 4.10)}$$

So $\log \frac{\mathfrak{m} \circ f}{\mathfrak{m}} \succ f - x \sim 1$. Since $f > x$, and $\mathfrak{m}$ is positive and infinitesimal, we have $\mathfrak{m} \circ f < \mathfrak{m}$. We deduce that $\log \frac{\mathfrak{m} \circ f}{\mathfrak{m}} < 0$ whence $\log \frac{\mathfrak{m} \circ f}{\mathfrak{m}} < \mathbf{R}$, so $\mathfrak{m} \circ f \prec \mathfrak{m}$. In particular, we have $\mathfrak{m} \circ f \in \mathbf{R}[\![\mathfrak{S}]\!]$. By Proposition 1.4, the family $\delta, \delta \circ f, \delta \circ f \circ f, \ldots$ is summable. The series

$$g := \delta + \delta \circ f + \delta \circ f \circ f + \cdots$$

satisfies $g \circ f + \delta = g$, whence $\Gamma(g) = g$ as desired. □

Proposition 4.7 follows from Propositions 4.18 and 4.19.

## 4.3 The case of finitely nested hyperseries

We now work within $\tilde{\mathbb{L}} = \mathbb{R}[\![\tilde{\mathfrak{L}}]\!]$. We write $\mathfrak{U}$ for the subgroup of $\tilde{\mathfrak{L}}$ generated by $\ell_0^{\mathbb{R}}$ and all groups $\widetilde{\mathfrak{L}_{<\omega^\nu}} \circ e_{\omega^\nu}^{\ell_0}$ for $\nu \in \mathbf{On}$ with $\nu > 0$, and set $\mathbb{U} := \mathbb{R}[\![\mathfrak{U}]\!]$. We have $\partial(\ell_0^{\mathbb{R}}) \subseteq \mathbb{R}\,\ell_0^{\mathbb{R}}$. Given an ordinal $\nu > 0$, the chain rule gives $(e_{\omega^\nu}^{\ell_0})' = \frac{1}{\ell'_{\omega^\nu}} \circ e_{\omega^\nu}^{\ell_0}$. Therefore

$$\begin{aligned}
\operatorname{supp} \partial(\widetilde{\mathfrak{L}_{<\omega^\nu}} \circ e_{\omega^\nu}^{\ell_0}) &\subseteq e_{\omega^\nu}^{\ell_0}{}' \partial(\widetilde{\mathfrak{L}_{<\omega^\nu}}) \circ e_{\omega^\nu}^{\ell_0} &&\text{(by the chain rule)} \\
&\subseteq (\widetilde{\mathfrak{L}_{<\omega^\nu}} \circ e_{\omega^\nu}^{\ell_0}) \cdot (\widetilde{\mathfrak{L}_{<\omega^\nu}} \circ e_{\omega^\nu}^{\ell_0}) \\
&\subseteq (\widetilde{\mathfrak{L}_{<\omega^\nu}}) \circ e_{\omega^\nu}^{\ell_0}.
\end{aligned}$$

We deduce with the Leibniz rule and by strong linearity of $\partial$ that $\partial(\mathbb{U}) \subseteq \mathbb{U}$.

**Lemma 4.20.** *For all $\nu_1, \ldots, \nu_p \in \mathbf{On}$ with $0 < \nu_1 < \cdots < \nu_p$ and all $(\mathfrak{m}_0, \ldots, \mathfrak{m}_n) \in \ell_0^{\mathbb{R}} \times \widetilde{\mathfrak{L}_{<\omega^\nu}} \circ e_{\omega^{\nu_1}}^{\ell_0} \times \cdots \times \widetilde{\mathfrak{L}_{<\omega^\nu}} \circ e_{\omega^{\nu_p}}^{\ell_0}$, we have $\mathfrak{m}_0 \cdots \mathfrak{m}_n \succcurlyeq \mathfrak{m}_i$ for $i = \max\{j \in \{0, \ldots, p\} : \mathfrak{m}_j \neq 1\}$.*

**Proof.** By Corollary 2.8(c), it is enough to note that $(\widetilde{\mathfrak{L}_{<\omega^\mu}})^{\neq} \circ e_{\omega^\mu}^{\ell_0} \prec\!\prec (\widetilde{\mathfrak{L}_{<\omega^\nu}})^{\neq} \circ e_{\omega^\nu}^{\ell_0}$ whenever $0 < \mu < \nu \in \mathbf{On}$ and that $\ell_0^{\mathbb{R}^\times} \prec\!\prec (\widetilde{\mathfrak{L}_{<\omega^\mu}})^{\neq} \circ e_{\omega^\mu}^{\ell_0}$. $\square$

**Corollary 4.21.** *For all $\mathfrak{u} \in \mathfrak{U}^{\neq}$, we have $\mathfrak{u} \succcurlyeq x$.*

By Lemma 3.17, the field $\mathbb{U}$ is a steep conjugate of $\tilde{\mathbb{L}}$. Theorem 4.6 applies and yields:

**Corollary 4.22.** *All positive elements in the group $\tilde{\mathbb{L}}^{>,\succ}$ are conjugate.*

**Remark 4.23.** In order to obtain a set-sized solution $\mathcal{G}$ to the conjugacy problem, it is enough to consider the closure under solutions $V$ of $V \circ f = V + 1$ for $f > \ell_0$, composition and inverses, of $\{\ell_0 + 1\}$. So set $\mathcal{G}_0 := \{\ell_0 + \mathbb{R}\}$, and for $n \in \mathbb{N}$, define $\mathcal{G}_{n+1}$ as the subgroup of $\tilde{\mathbb{L}}^{>,\succ}$ generated by $\mathcal{G}_n \cup \{V \in \tilde{\mathbb{L}}^{>,\succ} : \exists f \in \mathcal{G}_n, (f > \ell_0 \wedge V \circ f = V + 1)\}$. For any fixed $f > \ell_0$, the class of series $V$ with $V \circ f = V + 1$ is a set (see (4.3) below), so each $\mathcal{G}_n$ is a set. Thus $\mathcal{G} := \bigcup_{n \in \mathbb{N}} \mathcal{G}_n$ is a set-sized solution.

## 4.4 A simple inequality

Assume here that $\mathbb{S} = \mathbf{R}[\![\mathfrak{M}]\!]$, $\mathbb{T} = \mathbf{R}[\![\mathfrak{N}]\!]$ satisfy the properties listed in the beginning of Section 4 over an arbitrary ordered field of scalars $\mathbf{R}$, and that the group $\mathbb{S}^{>,\succ}$ solves the conjugacy problem. In other words, all $f, g > x$ are conjugate in $(\mathbb{S}^{>,\succ}, \circ, x)$. For instance we could take $\mathbb{S} = \tilde{\mathbb{L}}$. Our goal is to solve the inequality

$$f \circ g \geqslant f \circ g, \tag{4.2}$$

for $f, g \in \mathbb{S}^{>,\succ}$ with $f, g > x$. We start with the simpler case when $g = x + 1$.

**Lemma 4.24.** *Let $f \in \mathbb{S}^{>,\succ}$. We have $f \circ (x+1) > f + 1$ if $f > x + \mathbf{R}$, and $f \circ (x+1) < f + 1$ if $f < x + \mathbf{R}$.*

**Proof.** Suppose that $f > x + \mathbf{R}$. The inequality $f - x > \mathbf{R}$ implies by since $(\mathbb{S}, \partial)$ is an H-field that $f' > 1$. If $f' \sim 1$, then $f = x + \delta$ for $\delta := f - x > \mathbf{R}$. We have

$$f \circ (x+1) - f = x + 1 - x + \delta \circ (x+1) - \delta = 1 + \delta \circ (x+1) - \delta$$

By monotonicity, we have $\delta \circ (x+1) - \delta > 0$, whence $f \circ (x+1) > f + 1$. If $f' \asymp 1$ and $f' \not\sim 1$, then $f = r\,x + \varepsilon$ for a certain $r \in \mathbf{R}$ with $r > 1$ and a certain $\varepsilon \prec x$. By Lemma 4.4, we have $f \circ (x+1) - f \sim \tau_f \circ (x+1) - \tau_f \sim r$. We deduce that $f \circ (x+1) - f > 1$, hence the result.

Suppose now that $f' \succ 1$. By Lemma 4.4, we have $f \circ \left(x + \frac{2}{\tau'_f}\right) - f \sim \tau_f \circ \left(x + \frac{2}{\tau'_f}\right) - \tau_f$. Now $2/\tau'_f \prec \eth_f / \eth'_f$ so **TE** yields

$$\tau_f \circ \left(x + \frac{2}{\tau'_f}\right) - \tau_f \sim \tau'_f \frac{2}{\tau'_f} = 2.$$

In particular, $f \circ \left(x + \frac{2}{\tau'_f}\right) - f > 1$. Since $\tau'_f > \mathbf{R}$, we have $f \circ (x+1) - f > f \circ \left(x + 2/\tau'_f\right) - f$ by monotonicity, hence $f \circ (x+1) - f > 1$. The proof in the case when $f < x + \mathbf{R}$ follows from symmetric arguments. □

**Lemma 4.25.** *Let $f \in \mathbb{S}^{>,\succ}$ with $f \not< x + \mathbf{R}$ and $f \not> x + \mathbf{R}$. Then there is a unique $r_f \in \mathbf{R}$ such that $\delta := f - x + r_f$ is infinitesimal, and for all $r \in \mathbf{R}^>$, we have*

$$\begin{array}{rcll} f \circ (f + r) & < & f + r & \text{if } \delta > 0 \\ f \circ (f + r) & > & f + r & \text{if } \delta < 0. \end{array}$$

**Proof.** We must have $f - x \asymp 1$ so $f - x = r_f + \varepsilon$ for unique $r_f \in \mathbf{R}$ and $\varepsilon \prec 1$. Let $r \in \mathbf{R}^>$. So

$$f \circ (x + r) - f + g = r \circ (x + r) - \varepsilon.$$

If $\varepsilon > 0$, then $\varepsilon^{-1} \in \mathbb{S}^{>,\succ}$ so the function $\widehat{\varepsilon^{-1}}$ is strictly increasing on $\mathbb{S}^{>,\succ}$ by monotonicity. We deduce that $\varepsilon \circ (x+r) - \varepsilon < 0$, so $f \circ (x + r) < f + r$. The case when $\varepsilon < 0$ is symmetric. □

Recall that we write $\mathcal{C}(g) := \{h \in \mathbb{S}^{>,\succ} : h \circ g = g \circ h\}$ for the centralizer of a $g \in \mathbb{S}^{>,\succ}$.

**Corollary 4.26.** *We have $\mathcal{C}(x + 1) = x + \mathbf{R}$.*

**Proof.** We have $\mathcal{C}(x+1) \subseteq x + \mathbf{R}$ by Lemmas 4.24 and 4.25, whereas the converse inclusion $\mathcal{C}(x+1) \supseteq x + \mathbf{R}$ is immediate. □

## 4.5 Compositional iterates

Let $h \in \mathbb{S}^{>,\succ}$ with $h > x$ and let $U, V \in \mathbb{S}^{>,\succ}$ with

$$V \circ h = V + 1 \qquad \text{and} \qquad U \circ h = U + 1.$$

We have $V \circ h \circ V^{\mathrm{inv}} = x + 1 = U \circ h \circ U^{\mathrm{inv}}$, so $(U \circ V^{\mathrm{inv}}) \circ (x+1) = (x+1) \circ (U \circ V^{\mathrm{inv}})$. Thus by Corollary 4.26, there exists an $r_0 \in \mathbf{R}$ with

$$V = U + r_0. \tag{4.3}$$

So for all $r \in \mathbf{R}$, the series $V^{\mathrm{inv}} \circ (V + r)$ does not depend on the choice of $V$. We write

$$h^{[r]} := V^{\mathrm{inv}} \circ (V + r) \quad \text{for any } V \in \mathbb{S}^{>,\succ} \text{ with } V \circ h = V + 1.$$

We also write $(h^{\mathrm{inv}})^{[r]} := h^{[-r]}$ and $x^{[r]} = x$ for all $r \in \mathbf{R}$, so the operation

$$[\cdot] : \mathbf{R} \times \mathbb{S}^{>,\succ} \longrightarrow \mathbb{S}^{>,\succ} ; (r, h) \mapsto h^{[r]}$$

is well-defined. Forgetting the ordering for a moment, this endows $(\mathbb{S}^{>,\succ}, \circ, x)$ with a structure of **R**-group or **R**-exponential group as per Miasnikov and Remeslennikov [28]. In fact, in view of [27, Corollary 4.1], the ring **R** is the centroid of $(\mathbb{S}^{>,\succ}, \circ, x)$.

**Proposition 4.27.** *For all $h \in \mathbb{S}^{>,\succ} \setminus \{x\}$, the function $r \mapsto h^{[r]}$ is an isomorphism of ordered groups $(\mathbf{R}, +, 0, <) \longrightarrow (\mathcal{C}(h), \circ, x, <)$ with $h^{[1]} = h$.*

**Proof.** Fix a $V$ with $V \circ h = V + 1$. So $h = V^{\text{inv}} \circ (V + 1) = h^{[1]}$ by definition. For all $\varphi \in \mathbb{S}^{>,\succ}$ we have

$$\begin{aligned}\varphi \circ (x+1) = (x+1) \circ \varphi &\iff V^{\text{inv}} \circ \varphi \circ (V+1) = V^{\text{inv}} \circ (x+1) \circ \varphi \circ V \\ &\iff (V^{\text{inv}} \circ \varphi \circ V) \circ V^{\text{inv}} \circ (V+1) = V^{\text{inv}} \circ (V+1) \circ (V^{\text{inv}} \circ \varphi \circ V) \\ &\iff (V^{\text{inv}} \circ \varphi \circ V) \circ h = h \circ (V^{\text{inv}} \circ \varphi \circ V).\end{aligned}$$

Therefore $\mathcal{C}(h) = V^{\text{inv}} \circ \mathcal{C}(x+1) \circ V = V^{\text{inv}} \circ (x + \mathbf{R}) \circ V = \{h^{[r]} : r \in \mathbf{R}\}$. For all $r_0, r_1 \in \mathbf{R}$, we have $(x + r_0) \circ (x + r_1) = x + r_0 + r_1$ so $h^{[r_0 + r_1]} = h^{[r_0]} \circ h^{[r_1]}$. Furthermore, we have $r < r_1 \Longrightarrow V + r_0 < V + r_1 \Longrightarrow h^{[r_0]} < h^{[r_1]}$ by monotonicity. So $h \mapsto h^{[r]}$ is an isomorphism of ordered groups. □

In particular, centralizers of non-trivial elements in $\mathbb{S}^{>,\succ}$ are Abelian. For $n \in \mathbb{N}$, the series $h^{[n]}$ is the $n$-fold compositional iterate of $h$, whereas $h^{[-n]}$ is the $n$-fold compositional iterate of $h^{\text{inv}}$. Thus for $q \in \mathbb{Q}$, the series $h^{[q]}$ is a fractional iterate of $h$. For instance, we have a solution $h^{[1/2]}$ to the formal Schröder equation $y \circ y = h$. Proposition 4.28 below shows that $h^{[1/2]}$ is in fact its unique solution in $\mathbb{S}^{>,\succ}$.

**Proposition 4.28.** *For $h \in \mathbb{S}^{>,\succ}$ and $r_0, r_1 \in \mathbf{R}$. We have $(h^{[r_0]})^{[r_1]} = h^{[r_0 r_1]}$.*

**Proof.** We treat the case when $h > x$ and $r_0, r_1 > 0$. Let $V \in \mathbb{S}^{>,\succ}$ with $V \circ h = V + 1$ and set $U := (r_0^{-1} x) \circ V$. We have $(r_0^{-1} x) \circ (x + r_0) = (x + 1) \circ (r_0^{-1} x)$, so

$$\begin{aligned}U \circ h^{[r]} &= (r_0^{-1} x) \circ V \circ V^{\text{inv}} \circ (x + r_0) \circ V \\ &= (r_0^{-1} x) \circ (x + r_0) \circ V \\ &= U + 1.\end{aligned}$$

Therefore

$$\begin{aligned}(h^{[r_0]})^{[r_1]} &= U^{\text{inv}} \circ (U + r_1) \\ &= V^{\text{inv}} \circ (r_0 x) \circ (x + r_1) \circ (r_0^{-1} x) \circ V \\ &= V^{\text{inv}} \circ (x + r_0 r_1) \circ V \\ &= h^{[r_0 r_1]}.\end{aligned}$$

Since $f^{[-r]} = (f^{\text{inv}})^{[r]} = (f^{[r]})^{\text{inv}}$ for all $f \in \mathbb{S}^{>,\succ}$ and $r \in \mathbf{R}$, the other cases follow. □

We can now solve (4.2).

**Proposition 4.29.** *For $f, g > x$, the solution to the inequation $f \circ g \geqslant g \circ f$ is as follows:*

  a) *If $f > \mathcal{C}(g)$, then $f \circ g > g \circ f$.*
  b) *If $g > \mathcal{C}(f)$, then $f \circ g < g \circ f$.*
  c) *If $f \not> \mathcal{C}(g)$ and $g \not> \mathcal{C}(f)$, then for $g_f := \sup\{h \in \mathcal{C}(g) : h \leqslant f\}$, we have*

$$\begin{aligned}f \circ g &< g \circ f \quad \text{if } f > g_f \\ f \circ g &> g \circ f \quad \text{if } f < g_f, \text{ and} \\ f \circ g &= g \circ f \quad \text{if } f = g_f.\end{aligned}$$

**Proof.** Let $V \in \mathbb{S}^{>,\succ}$ with $V \circ g = V + 1$ and write $h := V \circ f \circ V^{\text{inv}}$. Note that $f \circ g \geqslant g \circ f \iff h \circ (x+1) \geqslant h+1$. So we may assume that $g = x+1$ and $f = h$. Then *a)* and *b)* follow from Lemma 4.24 and Proposition 4.27. Note that the element $r_f$ in Lemma 4.25 is the supremum of $\{r \in \mathbf{R} : x + r \leqslant f\}$. So *c)* follows from Lemma 4.25 and Proposition 4.27. □

**Proposition 4.30.** *Let $f, g \in \mathbb{S}^{>,\succ}$ with $f \geqslant g > x$. Then for all $g_0 \in \mathcal{C}(g)$, there is an $f_0 \in \mathcal{C}(f)$ with $f_0 \geqslant g_0$.*

**Proof.** We may assume that $f, g_0 > g$. We have $g_0 = g^{[r]}$ for a certain $r \in \mathbf{R}$ with $r \geqslant 1$, and we claim that $f^{[r]} > g_0$. Let $U, V \in \mathbb{S}^{>,\succ}$ with $U \circ f \circ U^{\mathrm{inv}} = x+1 = V \circ g \circ V^{\mathrm{inv}}$, so

$$f = U^{\mathrm{inv}} \circ (x+1) \circ U > V^{\mathrm{inv}} \circ (x+1) \circ V = g.$$

Writing $W := U \circ V^{\mathrm{inv}}$, we have $W^{\mathrm{inv}} \circ (x+1) \circ W > x+1$. Recall that

$$f^{[r]} = U^{\mathrm{inv}} \circ (x+r) \circ U \qquad \text{and} \qquad g_0 = g^{[r]} = V^{\mathrm{inv}} \circ (x+r) \circ V.$$

So $f^{[r]} > g_0$ if and only if $W^{\mathrm{inv}} \circ (x+r) \circ W > x+r$. Suppose that $W > x$. In view of Proposition 4.29, since $\mathcal{C}(x+r) = \mathcal{C}(x+1)$ and $r > 0$, we have $W^{\mathrm{inv}} \circ (x+r) \circ W > x+r$ if and only if $W^{\mathrm{inv}} \circ (x+1) \circ W > x+1$. If $W < x$, then apply Proposition 4.29 to $W^{\mathrm{inv}}$ instead. In any case, we get $f^{[r]} > g_0$ as claimed. □

**Acknowledgments.** We thank Gerald A. Edgar for his answers to our questions, Jean-Philippe Rolin and Tamara Servi for their advice and their help in reading parts of this paper, and Vincenzo Mantova for our discussions regarding compositional inverses.